\documentclass[reqno,10pt]{amsart}
\usepackage{amscd,amssymb,amsmath}
\usepackage{latexsym}
\usepackage[all]{xy}

\def\into{\hookrightarrow}

\def\toisom{\widetilde{\to}}

\def\.{,\dots ,}
\def\wt{\widetilde}
\def\wh{\widehat}
\def\ol{\overline}
\def\ul{\underline}
\def\wtimes{\wh{\otimes}}
\def\circcirc{{\circ\circ}}

\def\Mor{{\rm Mor}}

\def\Val{{\rm Val}}

\def\Val{{\rm Val}}
\def\Var{{\rm Var}}

\def\Spf{{\rm Spf}}
\def\Sp{{\rm Sp}}
\def\Gal{{\rm Gal}}

\def\Spa{{\rm Spa}}
\def\Coh{{\rm Coh}}

\def\Spec{{\rm Spec}}
\def\Proj{{\rm Proj}}
\def\Frac{{\rm Frac}}
\def\MSpec{\calM{\rm Spec}}
\def\MProj{\calM{\rm Proj}}
\def\bfSpec{{\bf Spec}}
\def\bfProj{{\bf Proj}}
\def\Bl{{\rm Bl}}

\def\hatBl{{\rm \wh{Bl}}}

\def\dim{{\rm dim}}

\def\rig{{\rm rig}}
\def\ad{{\rm ad}}

\def\an{{\rm an}}
\def\ar{{\rm ar}}

\def\Fsch{{\calF sch}}

\def\opp{{\rm op}}

\def\Int{{\rm Int}}
\def\Ker{{\rm Ker}}

\def\RZ{{\rm RZ}}
\def\Im{{\rm Im}}

\def\bir{{\rm bir}}

\def\cha{{\rm char}}

\def\rank{{\rm rank}}

\def\trdeg{{\rm tr.deg.}}

\def\bfA{{\bf A}}

\def\bfC{{\bf C}}

\def\bfF{{\bf F}}
\def\bfG{{\bf G}}

\def\bfN{{\bf N}}

\def\bfP{{\bf P}}
\def\bfQ{{\bf Q}}
\def\bfR{{\bf R}}

\def\bfZ{{\bf Z}}

\def\gtC{{\mathfrak C}}
\def\gtD{{\mathfrak D}}

\def\gtV{{\mathfrak V}}

\def\gtX{{\mathfrak X}}
\def\gtY{{\mathfrak Y}}
\def\gtZ{{\mathfrak Z}}

\def\calA{{\mathcal A}}
\def\calB{{\mathcal B}}
\def\calC{{\mathcal C}}

\def\calF{{\mathcal F}}
\def\calG{{\mathcal G}}
\def\calH{{\mathcal H}}

\def\calJ{{\mathcal J}}

\def\calM{{\mathcal M}}

\def\calO{{\mathcal O}}

\def\calX{{\mathcal X}}
\def\calY{{\mathcal Y}}

\def\oX{{\ol X}}

\def\ox{{\ol x}}

\def\uR{{\ul R}}
\def\uS{{\ul S}}
\def\uT{{\ul T}}

\def\ua{{\ul a}}

\def\uf{{\ul f}}
\def\ug{{\ul g}}

\def\ur{{\ul r}}
\def\us{{\ul s}}

\def\tilA{{\wt A}}

\def\tilS{{\wt S}}
\def\tilT{{\wt T}}

\def\tilX{{\wt X}}
\def\tilY{{\wt Y}}

\def\tilk{{\wt k}}
\def\till{{\wt l}}

\def\tilx{{\wt x}}
\def\tily{{\wt y}}

\def\hatA{{\wh A}}

\def\hatK{{\wh K}}

\def\hatX{{\wh X}}

\def\hatk{{\wh k}}

\def\tilcalA{{\wt\calA}}
\def\tilcalB{{\wt\calB}}

\def\hatcalA{{\wh\calA}}

\def\hatgtX{{\wh\gtX}}

\def\kcirc{k^\circ}
\def\lcirc{l^\circ}

\def\hatkcirc{\hatk^\circ}

\def\calAcirc{\calA^\circ}
\def\calBcirc{\calB^\circ}

\def\kcirccirc{k^{\circ\circ}}

\def\calAcirccirc{\calA^\circcirc}

\def\tilphi{{\wt\phi}}

\def\tilchi{{\wt\chi}}

\def\phicirc{\phi^\circ}

\def\chicirc{\chi^\circ}

\def\lam{{\lambda}}
\def\Lam{{\Lambda}}

\def\veps{\varepsilon}
\def\ve{\veps}

\def\wHx{{\widetilde{\calH(x)}}}

\def\wHz{{\widetilde{\calH(z)}}}

\def\whka{{\wh{k^a}}}

\def\R+*{{\bf R^*_+}}

\newtheorem{theorsect}{Theorem}[section]

\newtheorem{theor}{Theorem}[subsubsection]

\newtheorem{conj}[theor]{Conjecture}

\theoremstyle{definition}
\newtheorem{definsect}[theorsect]{Definition}

\newtheorem{defin}[theor]{Definition}
\newtheorem{rem}[theor]{Remark}
\newtheorem{fact}[theor]{Fact}
\newtheorem{exam}[theor]{Example}
\newtheorem{exer}[theor]{Exercise}
\newtheorem{exex}[theor]{Example/Exercise}
\newtheorem{defex}[theor]{Definition/Exercise}
\newtheorem{defexex}[theor]{Definition/Example/Exercise}

\begin{document}

\author{Michael Temkin}
\title{Introduction to Berkovich analytic spaces}
\address{\tiny{Einstein Institute of Mathematics, The Hebrew University of Jerusalem, Giv'at Ram, Jerusalem, 91904, Israel}}
\email{\scriptsize{temkin@math.huji.ac.il}}
\thanks{I want to thank A. Ducros for careful reading of the notes and making many valuable comments.}

\maketitle

\section{Introduction}
This paper presents an extended version of lecture notes for an introductory course on Berkovich analytic spaces that I gave in 2010 at Summer School "Berkovich spaces" at Institut de Mathématiques de Jussieu.

\subsection{Berkovich spaces and some history}

\subsubsection{Naive non-archimedean analytic spaces}
Since non-archimedean complete real valued fields (e.g. $\bfQ_p$) were discovered in the beginning of the last century, it was very natural to try to develop a theory of analytic spaces over such a field $k$ analogously to the theory of real or complex analytic spaces. At least, one can naturally define analytic functions on an open subset $V$ of the naive affine space $A^n_k=k^n$ as the convergent power series on $V$. This allowed to introduce some naive $k$-analytic spaces, but the theory was not rich enough. Actually, a global theory of such varieties does not make too much sense because the topology of $k$ is totaly disconnected. In particular, locally analytic (and even locally constant) functions do not have to be globally analytic.

\subsubsection{Rigid geometry}
In his study of elliptic curves with multiplicative bad reduction over a non-archimedean field $k$, Tate discovered about 1960 that these curves admit a natural uniformization by $G_m$. The latter was given as an abstract isomorphism of groups $k^\times/q^\bfZ\toisom E(k)$, and even such expert as Grothendieck doubted at first that this was not an accidental brute force isomorphism. Tate suspected, however, that his isomorphism can be interpreted as an analytic one, and he had to develop a good global theory of non-archimedean analytic spaces to make this
rigorous. This research resulted in Tate's definition of rigid geometry, whose starting idea was to simply forbid all bad open coverings (responsible for disconnectedness) and to shrink the set of analytic functions accordingly. As a result, one obtains a good theory of sheaves of analytic functions, but the underlying topological spaces have to be replaced with certain topologized (or Grothendieck) categories, also called $G$-topological spaces.

\subsubsection{Berkovich spaces}
More recently, some other approaches to non-archi\-me\-dean geometry were discovered: Raynaud's theory of formal models, Berkovich's analytic geometry and Huber's adic geometry. They all allow to define (nearly) the same categories of $k$-analytic spaces, but provide analogs of rigid spaces with additional structures invisible in the classical Tate's theory. Also, they extend the category of rigid spaces in different directions. Here we only discuss Berkovich's theory, which was developed in \cite{berbook} and \cite{berihes}. In this theory, classical rigid spaces are saturated with many new points (analogs of non-closed points of algebraic varieties), and the obtained spaces are honest topological spaces. In addition, the underlying topological spaces are rather nice (locally pathwise connected, for example).

Now, let us list some interesting features of Berkovich theory that distinguish it from all other approaches. First, one can work with all positive real numbers almost as well as with the values of $|k^\times|$. In particular, one can study rings of power series with radii of convergence linearly independent of $|k^\times|$. The latter fact allows to include the case of a trivially valued $k$ into the theory, and the theory of such $k$-analytic fields has already been applied to classical problems of algebraic geometry. Another interesting feature is that in a similar fashion one can develop (an equivalent form of) the usual theory of complex analytic spaces. Moreover, one can define Berkovich spaces that include both archimedean and non-archimedean worlds, for example, the affine line over $(\bfZ,|\ |_\infty)$.

\subsection{Structure of the notes}

We do not aim to prove all results discussed in these notes (and this is impossible in a six lecture long course). Our goal is to make the reader familiar with basic definitions, constructions, techniques and results of non-archimedean analytic geometry. Therefore, we prefer to formulate difficult results as Facts, and in some cases we discuss main ideas of their proofs in Remarks. Easy corollaries from these results (that may themselves be important pieces of the theory) are then suggested to the reader as exercises. Many exercises are provided with hints, but it may be worth to first try to solve them independently (especially, those not marked
with an asterisk).

\subsubsection{Overview}
The course is divided into five sections as follows. First, we study in \S\ref{2} semi-norms, norms and valuations, basic operations with these objects, Banach rings and their spectra. Then we describe the structure of $\calM(\bfZ)$, and after that we switch completely to the non-archimedean world. We finish the section with describing affine line over an algebraically closed non-archimedean field. In \S\ref{3} we introduce $k$-affinoid algebras and spaces and study their basic properties. In \S\ref{4}, this local theory is used to introduce and study global $k$-analytic spaces. Relations of $k$-analytic spaces with other categories are studied in \S\ref{5}. This includes analytification of algebraic varieties and GAGA, generic fibers of formal $\kcirc$-schemes and Raynaud's theory, and some discussion of rigid and adic geometries. In addition, in \S\ref{RZsec} we study local structure of analytic spaces by use of Riemann-Zariski (or birational) spaces over the residue field $\tilk$. Finally, in \S\ref{6} we study $k$-analytic curves in details. In particular, we describe their local and global structure and explain how this is related to the stable reduction theorem for formal $\kcirc$-curves.

\subsubsection{References and other sources}
The main references that helped me to prepare the course are \cite{berbook}, \cite{berihes} and \cite{bernotes}. The first two are a book and a large article in which the non-archimedean analytic spaces were introduced. The third one is a lecture notes of an analogous introductory course given by Berkovich in Trieste in 2009. I recommend the third source as an alternative (and shorter) expositional introductory text. It is worth to note that the first three sections of \cite{bernotes} and of these notes are parallel, but the exposition is (often but not always) rather different. Also, this text contains much more exercises and remarks, and this seriously
increases its length. Finally, the reader may wish to consult lecture notes \cite{con} on non-archimedean geometry (including the rigid geometry) by Brian Conrad.

\subsubsection{Novelties}
Some things in these notes are new. We develop a general theory of $H$-strict $k$-analytic spaces (or $k_H$-analytic) spaces that are built using radii of convergence from a group $|k^\times|\subseteq H\subseteq\bfR^\times_+$, with the extreme cases being strict and general $k$-analytic spaces. The main reference for $H$-strict theory is \cite{ct}. We give Berkovich's definition of $k$-analytic spaces that uses atlases but also show how one can define these spaces without atlases. This relies on the new Fact \ref{newfact}, which allows to characterize affinoid spaces and their morphisms as certain Banach ringed spaces and their morphisms. Another new result is Fact \ref{bddfact}, which asserts that for a non-trivially valued $k$ any $k$-homomorphism of $k$-affinoid algebras is bounded. Each new fact is followed by an exercise with a detailed hint on proving it.

\subsubsection{Conventions}
Throughout these notes {\em ring} always means a commutative ring with unity. For any field $k$ by $k^s$ and $k^a$ we denote its separable and algebraic closures, respectively. We will use underline to denote finite tuples of real numbers or of coordinates. For example, a polynomial ring $k[T_1\. T_n]$ will often be denoted as $k[\uT]$, where $\uT$ is the tuple $(T_1\. T_n)$ of coordinates. Also we will use the notation $\uT^i=T_1^{i_1}\dots T_n^{i_n}$ for $i\in\bfN^n$. For example, a power series $f(\uT)\in k[[\uT]]$ can be uniquely written as $\sum_{i\in\bfN^n}a_i\uT^i$ with $a_i\in k$.

\tableofcontents

\section{Norms, valuations and Banach rings}\label{2}

\subsection{Seminorms}

\subsubsection{Seminormed groups}

\begin{defex}\label{adddef}
(i) A {\em seminorm} on an abelian group $A$ is a function $|\ |\colon A\to\bfR_+$ which is {\em sub-additive}, i.e. $|a+b|\le|a|+|b|$, and satisfies $|0|=0$ and $|-a|=|a|$. A seminorm is a {\em norm} if its kernel is trivial. If the seminorm is fixed then we call $A$ a {\em seminormed group}. A seminorm is {\em non-archimedean} if it satisfies the strong triangle inequality $|a+b|\le\max(|a|,|b|)$.

(ii) The morphisms in the category of seminormed abelian groups are {\em bounded homomorphism}, i.e. homomorphisms $\phi{\colon}A\to B$ such that $\|\phi(a)\|\le C|a|$ for some fixed constant $C=C(\phi)$. In particular, $(A,|\ |)$ and $(A,\|\ \|)$ are isomorphic if and only if the seminorms $|\ |$ and $\|\ \|$ are {\em equivalent}, i.e. there exists a constant $C>0$ such that $|a|\le C\|a\|$ and $\|a\|\le C|a|$ for any $a\in A$.

(iii) Any quotient $A/H$ possesses a {\em residue seminorm} $\|\ \|$ given by $\|a+H\|=\inf_{h\in H}|a+h|$. A homomorphism of seminormed groups $\phi{\colon}A\to B$ is {\em admissible} if the residue seminorm on $\phi(A)$ is equivalent to the seminorm induced from $B$.

(iv) We provide a seminormed ring $A$ with the {\em semimetric} $d(a,b)=|a-b|$. The induced {\em seminorm topology} is the weakest topology for which the balls $B_{a,r}=\{x\in A|\ |x-a|< r\}$ are open. This topology distinguishes points (i.e. is $T_0$) if and only if the seminorm is a norm. Two seminorms are equivalent if and only if their induced topologies coincide. Any bounded homomorphism is continuous with respect to the seminorm topologies (see also Exercise \ref{contex}).

(v) The {\em separated completion} $\hatA$ of a seminormed group $A$ is the set of equivalence classes of Cauchy sequences in $A$. Use continuity to extend the group structure to $\hatA$ and show that $\hatA$ is a normed group, the natural map $A\to\hatA$ is an admissible homomorphism (called the {\em separated completion homomorphism}) and its kernel is $\Ker(|\ |)$. In particular, $A/\Ker(|\ |)$ is a normed group with respect to the residue seminorm.
\end{defex}

\begin{rem}
Usually, we will simply say "completion" in the sequel. Sometimes we will say "separated completion" in order to stress that the completion homomorphism may have a kernel.
\end{rem}

\subsubsection{Seminormed rings and modules}

\begin{defex}
(i) A {\em seminorm} (reps. {\em norm}) on a ring $A$ is a seminorm (resp. norm) on the additive group of $A$ which is {\em submultiplicative}, i.e. $|ab|\le|a||b|$. If $|\ |$ is multiplicative, i.e. $|ab|=|a||b|$ and $|1|=1$, then it is called a {\em real semivaluation} (resp. {\em real valuation}). If such a structure is fixed then the ring is called {\em seminormed}, {\em normed}, {\em real valued} or {\em real semivalued}, accordingly.

(ii) If $A$ is seminormed then a seminorm on an $A$-module $M$ is an additive seminorm $\|\ \|$ such that $\|am\|\le C|a|\|m\|$ for a fixed $C=C(M)$ and any $a\in A$ and $m\in M$.

(iii) Formulate and prove the analogs of all results/definitions from \S\ref{adddef}, including separated completions and admissible homomorphisms.
\end{defex}

\begin{rem}
(i) General non-archimedean semivaluations on rings are defined similarly but with values in $\{0\}\coprod\Gamma$, where $\Gamma$ is a totally ordered multiplicative abelian group. (Note that due to the strong triangle inequality, the definition makes sense even though there is no addition on $\{0\}\coprod\Gamma$.)

(ii) When studying general semivaluations one usually does not distinguish between the {\em equivalent} ones, i.e. semivaluations that admit an ordered isomorphism $i{\colon}\Im(|\ |)\toisom\Im(\|\ \|)$ such that $i\circ|\ |=\|\ \|$. This is the only reasonable possibility in the case when the group of values $\Gamma$ is not fixed. On the other side, it is very important that we do distinguish equivalent but not equal real semivaluations.

(iii) The valuation terminology is not unified in the literature. For example, in adic geometry of R.Huber, any semivaluation is called a valuation.
\end{rem}

In the following exercise we provide some definitions, examples and constructions related to seminormed rings.

\begin{defexex}\label{ostex}
Let $(A,|\ |)$ be a normed ring (the constructions make sense for seminormed rings but we will not need that).

(i) The {\em spectral seminorm} $\rho=\rho_\calA$ is the maximal {\em power-multiplicative} (i.e. $\rho(f^n)=\rho(f)^n)$) seminorm dominated by $|\ |$. Show that $\rho$ exists and is defined by $\rho(f)=\lim_{n\to\infty}|f^n|^{1/n}$.

(ii) For a tuple of positive numbers $\ur=(r_1\. r_n)$ provide $A[T_1\. T_n]$ with the norm
$$\|\sum_{i\in\bfN^n} a_i\uT^i\|^\ar_\ur=\sum_{i\in\bfN^n}|a_i|\ur^i$$ (where "ar" stands for archimedean) and let $A\{\ur^{-1}\uT\}^\ar=A\{r_1^{-1}T_1\. r_n^{-1}T_n\}^\ar$ denote its completion. This ring can be viewed as the ring of convergent power series over $A$ with polyradius of convergence $\ur$. Work this out: show that $A\{\ur^{-1}\uT\}^\ar$ is a subring of $\hatA[[\uT]]$ defined by a natural convergence conditions.

(iii) If $(M,|\ |_M)$ and $(N,|\ |_N)$ are normed $A$-modules (resp. rings) then we provide $M\otimes_A N$ with the {\em tensor product seminorm} $\|x\|=\inf(\sum_{i=1}^n|m_i|_M|n_i|_N)$ where the infimum is taken over all representations of $x$ of the form $x=\sum_{i=1}^n m_i\otimes n_i$. The separated completion of this seminormed module is denoted $M\wtimes^{\ar}_A N$ and called the (archimedean) {\em completed tensor product} of modules (resp. rings). We will later see that the tensor product seminorm is often not a norm.

(iv) The {\em trivial semi-norm} $|\ |_0$ on a ring $A$ sends $A\setminus\{0\}$ to $1$. It is power-multiplicative (resp. a valuation) if and only if $A$ is reduced (resp. integral).

(v) For any natural $n>1$ define the $n$-adic norm on $\bfQ$ by $|x|_n=n^d$, where $d\in\bfZ$ is the minimal number with $xn^d\in\bfZ_{(n)}$ (the localization of $\bfZ$ by all primes coprime with $n$). This norm is a valuation only when $n$ is prime. The equivalence class of $|\ |_n$ depends only on the set $p_1\. p_m$ of prime divisors of $n$. The completion $\bfQ_n$ with respect to $|\ |_n$ is called the ring of $n$-adic numbers. Show that $\bfQ_n=\oplus_{i=1}^m\bfQ_{p_i}$. In particular, the completion operation does not preserve the property of being an integral domain. Show that $\bfQ_p$ is a field. It is called the {\em field of $p$-adic numbers}.

(vi) Define $t$-adic valuations on $k[t]$ analogously to the $p$-adic valuation (they are trivial on $k$ and are uniquely determined by $r=|t|\in(0,1)$). Show that $k[[t]]$ is the completion.

(vii) Ostrowski's theorem provides a complete list of real semivaluations on $\bfZ$: the trivial valuation $|\ |_0$, the $p$-adic valuations $|\ |_{p,r}=(|\ |_p)^r$ for any $r\in(0,\infty)$, the archimedean valuations $|\ |_{\infty,r}=(|\ |_\infty)^r$ for $r\in(0,1]$ (where $|x|_\infty$ is the usual absolute value of $x$), and the semivaluations $|\ |_{p,\infty}$ that take $p\bfZ$ to $0$ and everything else to $1$.
\end{defexex}

\begin{rem}\label{anrem}
There is a certain analogy, that will be used later, between ideals on rings and bounded semi-norms on semi-normed rings. Exercises (iv) and (v) above indicate that multiplicative (resp. power-multiplicative) semi-norms correspond to prime (resp. reduced) ideals. In the style of the same analogy, passing from a semi-normed ring $(A,|\ |)$ to $(A/\Ker(\rho_A),\rho_A)$ can be viewed as an analog of reducing the ring (i.e. factoring a ring by its radical). The elements in the kernel of $\rho_A$ are called {\em quasi-nilpotent}
elements.
\end{rem}

\subsection{Banach rings and their spectra}
\subsubsection{Banach rings, algebras and modules}
\begin{defin}
(i) A {\em Banach ring} is a complete normed ring $\calA$ (i.e. the completion homomorphism $\calA\to\hatcalA$ is an isomorphism). A {\em Banach $\calA$-algebra} is a Banach algebra $\calB$ with a bounded homomorphism $\calA\to\calB$.

(ii) A {\em Banach $\calA$-module} is a complete normed $\calA$-module.
\end{defin}

Instead of the polynomial rings and tensor products of modules, when working with Banach rings and modules we will use the convergent power series rings and completed tensor products.

\begin{fact}
The valuation of any complete real valued field $k$ uniquely extends to any algebraic extension of $k$.
\end{fact}

\begin{exex}\label{contex}
Let $k$ be a complete real valued field and let $f{\colon}\calA\to\calB$ be a homomorphism of Banach $k$-modules (also called Banach $k$-spaces).

(i) Assume that $k$ is not trivially valued (e.g. $\bfR$, $\bfC$, $\bfQ_p$ or $\bfC((t))$). Then $f$ is bounded if and only if it is continuous.

(ii) Assume that $k$ is trivially valued. Show that as an abstract ring, $k\{r^{-1}T\}^\ar$ is isomorphic to $k[T]$ when $r\ge 1$ and is isomorphic to $k[[T]]$ when $r<1$. In particular, there exist continuous but not bounded homomorphisms of Banach $k$-algebras.
\end{exex}

\subsubsection{The spectrum}
The analogy from Remark \ref{anrem} suggests the following definition.

\begin{defin}
(i) {\em Spectrum} of a Banach ring $\calA$ is the set $\calM(\calA)$ of all bounded real semivaluations $|\ |_x$ on $\calA$ (i.e. $|\ |_x\le C|\ |$ for some $C$) provided with the weakest topology making continuous the maps $|f|{\colon}\calM(\calA)\to\bfR_+$ for all $f\in\calA$. The latter maps take $|\ |_x$ to $|f|_x$ and usually we will use the notation $x\in\calM(\calA)$ and $|f(x)|$ instead of $|\ |_x$ and $|f|_x$.

(ii) For any point $x\in\calM(\calA)$ the kernel of $|\ |_x$ is a prime ideal and hence $A/\Ker(|\ |_x)$ is an integral valued ring. The completed fraction field of this ring is called the {\em completed residue field of $x$} and we denote it as $\calH(x)$. The bounded character corresponding to $x$ will be denoted $\chi_x{\colon}\calA\to\calH(x)$.
\end{defin}

The following exercise shows that the definition of $\calM$ is analogous to the definition of $\Spec$.

\begin{exer}
The points of $\calM(\calA)$ are the isomorphism classes of bounded homomorphisms $\chi{\colon}\calA\to k$ whose image is a complete real valued field generated by the image of $\chi$ (i.e. $\Im(\chi)$ is not contained in a proper complete subfield of $k$).
\end{exer}

Here are basic facts about the spectrum. As one might expect, the general line of the proof is to construct enough points by use of Zorn's lemma.

\begin{fact}
(i) Let $\calA$ be a Banach ring. The spectrum $X=\calM(\calA)$ is compact, and it is empty if and only if
$\calA=0$.

(ii) The maximum modulus principle: $\rho(f)=\max_{x\in X}|f(x)|$.
\end{fact}

\begin{exer}\label{calmex}
(i) Extend $\calM$ to a functor to topological spaces, that is, for any bounded homomorphism of Banach rings $\phi{\colon}\calA\to\calB$ construct a natural continuous map $\calM(\phi){\colon}\calM(\calB)\to\calM(\calA)$.

(ii) An element $f\in\calA$ is invertible if and only if $\inf_{x\in X}|f(x)|>0$.

(iii) If $\calA$ is a Banach $k$-ring for a complete real valued field $k$ then $$\calM(\calA\wtimes_k k^a)/\Gal(k^s/k)\toisom\calM(A)$$ where $k^a$ and $k^s$ are provided with the extended valuation.

(iv)* Let $\calA$ be finite over $\bfZ$. Show that up to equivalence there exists unique structure of a Banach $(\bfZ,|\ |_\infty)$-algebra on $\calA$ and describe $\calM(\calA)$ similarly to the description of $\calM(\bfZ)$ in Exercise \ref{ostex}(vii). Analyze similarly the spectra of finite Banach $(k[T],|\ |_0)$-algebras, where $k$ is a field and $|\ |_0$ is the trivial valuation. (Hint: take the scheme $\Spec(\calA)$. Keep all its closed points, and for any generic point $x$ of a curve component of $\Spec(\calA)$ replace $x$ with all valuations on $k(x)$.)
\end{exer}

\subsubsection{Relative affine spectrum}
The following definition is not standard, but it seems to be convenient.

\begin{defin}
Let $\calA$ be a Banach ring and let $C$ be an $\calA$-algebra (without any norm). We define the analytic spectrum of $C$ as the set $\MSpec(C)$ of all real semivaluations on $C$ that are bounded on $\calA$ (i.e. the restriction of $|\ |_x$ to $\calA$ is bounded). Naturally, $\MSpec(C)$ is provided with the weakest topology making continuous each map $x\mapsto |f(x)|$ with $f\in\calA$.
\end{defin}

\begin{rem}
There is a natural projection $\MSpec(C)\to\calM(\calA)$, which is typically a non-compact map. We stress that $\MSpec(C)$ depends on the structure of $C$ as an $\calA$-algebra. In some sense, it is an analog of the relative $\bfSpec$ construction in algebraic geometry.
\end{rem}

\subsubsection{The affine space}
\begin{defin}
The {\em $n$-dimensional affine space} over a Banach ring $\calA$ is the topological space $\MSpec(\calA[T_1\. T_n])$.
\end{defin}

\begin{exer}
The affine space $\bfA^n_\calA$ is the union of closed {\em $\calA$-polydiscs} $\calM(\calA\{\ur^{-1}\uT\}^\ar)$ of {\em polyradius} $\ur=(r_1\. r_n)$. In particular, it is locally compact.
\end{exer}

\begin{rem}
For a complete real valued field $k$ one can provide $\bfA^n_k$ with the sheaf of analytic functions which are local limits of rational functions from $k(\uT)$. If $U\subset\bfA^n_k$ is open and $V$ is a Zariski closed subset of $U$ given by vanishing of analytic functions $f_1\. f_n$, then factoring by the ideal generated by $f_i$ one obtains a sheaf of analytic functions on $V$. Gluing such local models $V$ with the sheaves of analytic functions one can construct a theory of $k$-analytic spaces without boundary. The advantage of this approach is that it works equally well over $\bfC$ and over $\bfQ_p$. The main disadvantage of this approach is that it does not treat well enough the cases of a trivially valued $k$ and of analytic $\bfQ_p$-spaces with boundary. More details about the outlined approach can be found in \cite[\S1.5]{berbook} and \cite[\S1.3]{bernotes}. We will use another approach to construct non-archimedean analytic spaces.
\end{rem}

\subsection{Non-archimedean setting}

\subsubsection{Strong triangle inequality}
\begin{defin}
(i) A {non-archimedean seminorm} (resp. norm, semivaluation, etc.) is a seminorm $|\ |$ that satisfies the {\em strong triangle inequality} $|a+b|\le\max(|a|,|b|)$.

(ii) A {\em non-archimedean field} is a complete real valued field $k$ whose valuation is non-archimedean.
\end{defin}

\begin{exam}
(i) K\"ursch\'ak's proved that any complete real valued field containing $\bfC$ coincides with it (Gel'fand-Mazur proved a stronger claim later). Now, applying Ostrowski's classification we obtain that any complete real valued field, excluding $\bfR$ and $\bfC$, is non-archimedean.

(ii) For any ring $A$ its trivial seminorm is non-archimedean.
\end{exam}

\begin{exer}
Show that for an archimedean $k$ one has that $\bfA^1_k\toisom k^a/\Gal(k^a/k)$ where the image is provided with the valuation topology, i.e. $\bfA^1_k$ coincides with the naive affine line.
\end{exer}

In the sequel we will work only with non-archimedean seminorms, semivaluations, etc., so the word "non-archimedean" will usually be omitted. Let $\calA$ be a non-archimedean Banach ring. The basic definitions should now be adjusted as follows.

\begin{defexex}\label{exx}
(i) Check that $\calM(\calA)$ is the set of all non-archimedean bounded semivaluations on $\calA$.

(ii) The spectral seminorm $\rho_\calA$ is non-archimedean.

(iii) Non-archimedean definitions of $\calA\{\ur^{-1}\uT\}$, $M\wtimes_\calA N$ and their norms copy their archimedean analogs with maxima used instead of sums. For example, $$\|\sum_{i\in\bfN^n} a_i\uT^i\|_\ur=\max_{i\in\bfN^n} |a_i|\ur^i$$ Check that $\|\ \|_\ur$ is a valuation.

(iv) Calculus student's dream: a sequence $a_n$ in $\calA$ is Cauchy if and only if
$\lim_{n\to\infty}|a_n-a_{n+1}|=0$. In particular, a series $\sum_{n=0}^\infty a_n$ converges in $\calA$ if and only if $\lim_{n\to\infty}|a_n|=0$.
\end{defexex}

At this point we fix a non-archimedean ground field $k$ and start to develop non-archimedean analytic geometry over $k$. When developing this theory we will compare it from time to time to the classical theory of algebraic varieties over a field. Analogously to the latter theory, we will first introduce $k$-Banach algebras of topologically finite type and their spectra, called $k$-affinoid algebras and spaces. Then we will construct general $k$-analytic spaces by pasting $k$-affinoid ones. Despite this general similarity, many details in our theory are subtler. We will try to indicate critical moments where the theories differ.

\subsubsection{Reduction ring}
\begin{defin}
It follows from Exercise \ref{exx} that any non-archimedean ring $\calA$ contains an open subring $\calAcirc=\{a\in\calA|\ \rho(a)\le 1\}$ with an ideal $\calAcirccirc:=\{a\in\calA|\ \rho(a)<1\}$. The ring $\tilcalA=\calAcirc/\calAcirccirc$ is called the {\em reduction ring} of $\calA$.
\end{defin}

\begin{exex}
If $k$ is a non-archimedean field then $\kcirc$ is its valuation ring with maximal ideal $\kcirccirc$ and residue field $\tilk$.
\end{exex}

\subsubsection{Description of points of $\bfA^1_k$}\label{aflinesec}
\begin{defin}
Let $l$ be a non-archimedean $k$-field. Recall that $e_{l/k}$ is the cardinality of $|l^\times|/|k^\times|$ (may be infinite) and $f_{l/k}=[\till:\tilk]$. Transcendental analogs of these cardinals are $E_{l/k}=\rank_\bfQ(|l^\times|/|k^\times|\otimes_\bfZ\bfQ)$ and $F_{l/k}=\trdeg(\till/\tilk)$.
\end{defin}

\begin{exer}
Assume that $l$ is algebraic over the completion of its subfield $l_0$ which is of transcendence degree $n$ over $k$. Prove Abhyankar inequality: $E_{l/k}+F_{l/k}\le n$.
\end{exer}

We will use this result to classify points on $\bfA^1_k=\MSpec(k[T])$ (and a similar argument classifies points on any $k$-analytic curve).

\begin{defex}\label{typesex}
(0) A point $x\in\bfA^1_k$ is {\em Zariski closed} if $|\ |_x$ has a non-trivial kernel. Show that this happens if and only if $\calH(x)$ is finite over $k$. Show that otherwise $\calH(x)$ is a completion of $k(T)$ and to give a point which is not Zariski closed is the same as to give a real valuation on $k(T)$ that extends that of $k$. In particular, $E_{\calH(x)/k}+F_{\calH(x)/k}\le 1$.

(1) $x$ is of type 1 if $\calH(x)\subseteq\whka$.

(2) $x$ is of type 2 if $F_{\calH(x)/k}=1$.

(3) $x$ is of type 3 if $E_{\calH(x)/k}=1$.

(4) $x$ is of type 4 if $E_{\calH(x)/k}=F_{\calH(x)/k}=0$ and $x$ is not of type 1.

(5)* Show that in case (1) $\calH(x)$ may contain an infinite algebraic extension of $k$ (e.g. if $k=\bfQ_p$ then it may coincide with $\bfC_p=\wh\bfQ_p^a$). In particular, the map $\bfA^1_{\whka}\to\bfA^1_k$ usually has infinite (pro-finite) fibers. (Hint: fix elements $x_i\in k^s$ and take $T=\sum_{i=1}^\infty a_ix_i$ where $a_i\in k$ converge to zero fast enough; then use Krasner's lemma to show that $k(x_i)\subset\wh{k(T)}$.)
\end{defex}

\begin{rem}
More generally, $\calH(x)$ may contain an infinite algebraic extension of $k$ for type 4 points, but not for type 2 or 3 points.
\end{rem}

Assume now that $k$ is algebraically closed and let us describe the points of $\bfA^1_k$ in more details. By closed disc

\begin{exer}
(i) By {\em closed disc} $E(a,r)\subset \bfA^1_k$ of radius $r$ and with center at $a$ we mean the set of points of $\bfA^1_k$ that satisfy $|(T-a)(x)|\le r$. Show that $E(a,r)=\calM(k\{r^{-1}(T-a)\})$

(ii) Show that a valuation on $k[T]$ is determined by its values on the elements $T-a$ with $a\in k$. The number $r=\inf_{a\in k}|T-a|$ is called the {\em radius} of $x$ (with respect to the fixed coordinate $T$).

(iii) Assume that the infimum $r$ is achieved, say $r=|T-a|$. Show that

(a) if $r=0$ then $x$ is Zariski closed and of type 1 and $|f(x)|=|f(a)|$ for any $f\in k[T]$.

(b) if $r>0$ then $x$ is the maximal point of the disc $E(a,r)$ (i.e. $|f(x)|\ge|f(y)|$ for any $y\in E(a,r)$ and $f=\sum_{i=0}^nf_iT^i\in k[T]$) and $|f(x)|=\max_i |f_i|r^i$. If $r$ is {\em rational} in the sense that $r^n\in|k^\times|$ for some integral $n>0$ then $x$ is of type $2$, and otherwise $x$ is of type 3.

(iii) Assume that the infimum is not achieved, say $a_j\in k$ are such that the sequence $r_j=|T-a_j|$ decreases and tends to $r$. Then $x$ is of type 4, it is the only point in the intersection of the discs $E(a_j,r_j)$, and $|f(x)|=\inf_j|f(x_j)|$ where $x_j$ is the maximal point of $E(a_j,r_j)$. In particular, for an algebraically closed ground field $k$, type 4 points exist in $\bfA^1_k$ if and only if $k$ is not {\em spherically complete}, i.e. there exist nested sequences of discs over $k$ without common $k$-points.
\end{exer}

Actually, $\bfA^1_k$ is a sort of an infinite tree whose leaves are type 1 and 4 points.

\begin{exer}
(i) Use the previous exercise to prove that $\bfA^1_k$ is pathwise connected and simply connected. Moreover, show that for any pair of points $x,y\in\bfA^1_k$ there exists a unique path $[x,y]$ that connects them. (Hint: $[x,y]=[x,z]\cup[z,y]$ where $z$ is the maximal point of the minimal disc containing both $x$ and $y$ and the open path $(x,z)$ (resp. $(z,y)$) consists of the maximal points of discs that contain $x$ but not $y$ (resp. $y$ but not $x$).)

(ii) Show that $\bfA^1_k\setminus\{x\}$ is connected whenever $x$ is of type 1 or 4, consists of two components when $x$ is of type 3, and consists of infinitely many components naturally parameterized by $\bfP^1_\tilk$ when $x$ is of type 2. Thus, $\bfA^1_k$ is an infinite tree with infinite ramification at type 2 points. If $k$ is trivially valued then there is just one type 2 point and no type 4 points, so the tree looks like a star whose rays connect the type 2 point (the trivial valuation) with the Zariski closed points.
\end{exer}

Almost all non-discretely valued fields are not spherically complete.

\begin{exer}
(i) Let $k_0$ be trivially valued and let $k$ be the $t$-adic field $k_0((t))$. Show that $\bfC_p$ and $\whka$ are not spherically complete. (Hint: for example, choose centers of the discs at $\sum_{i=0}^n t^{l_i}$, where $l_i$ is a decreasing sequence of rational numbers that tends to a positive number $r$.)

(ii) Prove by Zorn's lemma that spherically complete, algebraically closed, and non-trivially valued non-archimedean fields exist.

(iii) Here is the only known explicit construction of such fields. Let $k_0$ be an algebraically closed trivially valued field and let $\Gamma\subseteq\bfR^\times_+$ be a divisible subgroup. Let $k=k_0((t^\Gamma))$ be the set of all series $\sum_{\gamma\in\Gamma}a_\gamma t^\gamma$ where $a_\gamma\in k_0$ and any increasing family of $\gamma$'s with non-zero $a_\gamma$ is finite (finite sums form the group ring $k_0[t^\Gamma]$). Show that one can naturally define multiplication that makes $k$ to a spherically complete and algebraically closed non-archimedean field with group of values $\Gamma$ and residue field $k_0$.
\end{exer}

\begin{rem}
The construction from (iii) is very nice, but I do not know about any application of the fields $k_0((t^\Gamma))$ to non-archimedean geometry. However, existence of a spherically complete closure plays important role in non-archimedean geometry. For example, few approaches to stable reduction theorem first prove the result over a spherically complete field, thus avoiding some troubles caused by type 4 points, and then establish the general case by a descent argument. It seems that the first such proof is due to van der Put. A similar strategy is also used in the recent work \cite{HL} by Hrushovski-Loeser, that we will recall in \S\ref{topsec}.
\end{rem}

\section{Affinoid algebras and spaces}\label{3}

\subsection{Affinoid algebras}

\subsubsection{The definition}
\begin{defin}\label{affdef}
(i) A {\em $k$-affinoid algebra} $\calA$ is a Banach $k$-algebra that admits an admissible surjective homomorphism from a Banach algebra of the form $k\{\ur^{-1}\uT\}$. We say that $\calA$ is {\em strictly $k$-affinoid} if one can choose $r_i\in|k^\times|$. More generally, we say that $\calA$ is {\em $H$-strict} for a group $|k^\times|\subseteq H\subseteq\bfR^\times_+$ if one can choose such a homomorphism with $r_i\in H$.

(ii) The category of (resp. $H$-strict, resp. strictly) $k$-affinoid algebras with bounded morphisms is denoted $k$-$AfAl$ (resp. $k_H$-$AfAl$, resp. $st$-$k$-$AfAl$). It will also be convenient to say {\em $k_H$-affinoid algebra} instead of $H$-strict $k$-affinoid algebra.
\end{defin}

\begin{exer}
Check that $H$-strictness depends only on the group $\sqrt{H}$ consisting of all elements $h^{1/n}$ with $h\in H$ and integral $n\ge 1$.
\end{exer}

\begin{rem}
The group $\sqrt{H}$ is not dense in $\bfR_+^\times$ if and only if $H=1$, and $1$-strict spaces are precisely the strictly analytic spaces over a trivially valued field. The case of $H=1$ is degenerate and often demonstrates a very special behavior. We will ignore it in all cases when it requires a separate argument.
\end{rem}

\begin{exex}\label{Krex}
Let $\ur=(r_1\. r_n)$ be a tuple of positive real numbers linearly independent over $|k^\times|$. Show that the $k$-affinoid ring $$K_\ur:=k\{\ur^{-1}\uT,\ur\ \uT^{-1}\}=k\{\ur^{-1}\uT,\ur\ \uS\}/(T_1S_1-1\. T_nS_n-1)$$ is a field and $K_\ur\toisom K_{r_1}\wtimes_kK_{r_2}\wtimes_k\dots\wtimes_kK_{r_n}$.
\end{exex}

\subsubsection{Basic properties}
Here is a summary of basic properties of $k$-affinoid algebras. Excellence was proved very recently by Ducros in \cite{duc2} (and the strictly affinoid case is due to Kiehl).

\begin{fact}\label{basicaffact}
(i) Any affinoid algebra $\calA$ is noetherian, excellent and all its ideal are closed.

(ii) If $f\in\calA$ is not nilpotent then there exists $C>0$ such that $\|f^n\|\le C\rho(f)^n$ for all $n\ge 1$. In particular, $f$ is not quasi-nilpotent (i.e. $\rho(f)>0$), and so $\rho$ is a norm if and only if $\calA$ is reduced.

(iii) If $\calA$ is reduced then the Banach norm on $\calA$ is equivalent to the spectral norm.

(iv) $\calA$ is $H$-strict if and only if $\rho(\calA)\subseteq\{0\}\cup\sqrt H$.
\end{fact}

In particular, $(A/\Ker(\rho_A),\rho_A)$ is equivalent to the quotient of $A$ by its radical (provided with the residue semi-norm). In view of Remark \ref{anrem}, this can be interpreted as equivalence of the "topological reduction" of $A$ and the usual reduction of $A$. The following example shows that even naively looking $k$-Banach algebras do not have to satisfy the same nice conditions.

\begin{exex}
Let $k$ be complete non-perfect field with a non-trivial valuation (e.g. $k=\bfF_p((t))$). Take any element $x$ lying in the completion of the perfect closure of $k$ and non-algebraic over $k$ (e.g. $x=t^{1+1/p}+t^{2+2/p^2}+\dots$) and let $K$ be the closure of $k(x)$ in $\whka$. (Note that $K=\calH(z)$ for a non Zariski closed point $z\in\bfA_k^1$ of type 1.) Show that the element $1\otimes x-x\otimes 1$ is a quasi-nilpotent element of $K\wtimes_k K$ which is not nilpotent.
\end{exex}

For strictly affinoid algebras one can say more.

\begin{fact}\label{bfact}
Let $\calA$ be a strictly $k$-affinoid algebra (where the trivially valued case is allowed).

(i) Noether normalization: there exists a finite admissible injective homomorphism $k\{T_1\. T_n\}\to\calA$.

(ii) Hilbert Nullstellensatz: $\calA\neq 0$ has a point in a finite extension of $k$.

(iii) The rings $k\{T_1\. T_n\}$ are of dimension $n$.
\end{fact}

\begin{rem}
(i) A very systematic and detailed theory of strictly affinoid algebras is developed in chapters 5 and 6 of \cite{BGR} (they are called "affinoid algebras" in loc.cit.). In particular, loc.cit. contains the proof of all claims of Facts \ref{basicaffact} and \ref{bfact} excluding the excellence result. Summarizing in a couple of words, first one develops a Weierstrass theory (preparation and division theorems) for strictly affinoid algebras. As a corollary one deduces analogs of two famous theorems about affine algebras: Noether normalization and Hilbert Nullstellensatz.  All these results are used to establish Fact \ref{basicaffact} in the strict case.

(ii) Berkovich introduced non-strict algebras in \cite{berbook} and suggested the following descent trick to deal with them. Obviously, for any $k$-affinoid algebra $\calA$ its base change $\calA\wtimes_k K_\ur$ is strictly $K_\ur$-affinoid for an appropriate $K_\ur$. This allows to show that many good properties known to hold for strictly $K_\ur$-affinoid algebras also hold for general $k$-affinoid algebras. In particular, this approach provides a simple reduction of all claims of Fact \ref{basicaffact}, excluding excellence, to the known strictly affinoid case.

(iii) It was not studied in the literature whether one can develop the whole theory for all affinoid algebras. My expectations are as follows. Weierstrass theory can be developed for all affinoid algebras. Hilbert Nullstellensatz holds in a corrected form that any affinoid $\calA$ has a point in a finite extension of some $K_r$, see \cite[Th. 2.7]{duc1}. I expect that the following weak form of Noether normalization is the best one can get (see Example \ref{irrex}(ii)). There exists injective homomorphisms $f{\colon}k\{r_1^{-1}T_1\. r_n^{-1}T_n\}\to\calA'$ and $g{\colon}\calA'\to\calA$ such that $f$ is finite admissible and $g$ has dense image (then $\calM(\calA)$ is a Weierstrass domain in a finite surjective covering $\calM(\calA')$ of a polydisc, as we will later see).
\end{rem}

Fact \ref{basicaffact} has the following corollary, which is very important when developing the theory of affinoid spaces.

\begin{exer}
Assume that $\phi{\colon}\calA\to\calB$ is a bounded homomorphism of $k$-affinoid algebras, $f_1\. f_n\in\calB$ are elements and $r_1\. r_n>0$ are real numbers. Then $\phi$ extends to a bounded homomorphism $\psi{\colon}\calA\{r_1^{-1}T_1\. r_n^{-1}T_n\}\to\calB$ with $\psi(T_i)=f_i$ if and only if $\rho_\calB(f_i)\le r_i$.
\end{exer}

\begin{defin}
Let $\calA$ be a $k$-affinoid algebra. Any Banach $\calA$-algebra that admits an admissible surjective homomorphism from $\calA\{r_1^{-1}T_1\. r_n^{-1}T_n\}$ is called {\em $\calA$-affinoid}. Obviously, it is also a $k$-affinoid algebra.
\end{defin}

\subsubsection{Finite $\calA$-modules}
It turns out that the theory of finite Banach $\calA$-modules is essentially equivalent to the theory of finite $\calA$-modules.

\begin{defin}
A Banach $\calA$-module $M$ is {\em finite} if it admits an admissible surjective homomorphism from a free module $\calA^n$ provided with the norm $||(a_1\. a_n)||=\max_{1\le i\le n}|a_i|$.
\end{defin}

\begin{fact}\label{basicmodfact}
(i) The categories of finite Banach $\calA$ modules and finite $\calA$-modules are equivalent via the forgetful functor. In particular, any $\calA$-linear map between finite Banach $\calA$-modules is admissible.

(ii) Completed tensor product with a finite Banach $\calA$-module $M$ coincides with the usual tensor product. Namely, $M\otimes_\calA N\toisom M\wtimes_\calA N$ for any Banach $\calA$-module $N$.
\end{fact}

\begin{exer}
Formulate and prove an analog of Fact \ref{basicmodfact} for the category of finite $\calA$-algebras. In addition, prove that any finite Banach $\calA$-algebra is $\calA$-affinoid.
\end{exer}

\subsubsection{Complements}

\begin{fact}
(i) Fibred coproducts exist in the category $k_H$-$AfAl$ and coincide with completed tensor products.

(ii) For any non-archimedean $k$-field $K$, the correspondence $\calA\mapsto\calA\wtimes_k K$ provides a ground field extension functor $k_H$-$AfAl\to K_{H|K^\times|}$-$AfAl$ compatible with completed tensor products.
\end{fact}

\begin{fact}\label{bddfact}
The ground field $k$ is not trivially valued if and only if any homomorphism between $k$-affinoid algebras is bounded.
\end{fact}

The latter fact was known only for strictly affinoid algebras, so we suggest a proof below.

\begin{exer}
(i) Show that any automorphism of the $k$-field $K_r$ from Example \ref{Krex} is bounded if and only if $k$ is not trivially valued. (Hint: you have to use arithmetical properties of $K_r$ because $\wh{K_r^a}$ obviously has a lot of non-bounded automorphisms.)

(ii)* Prove Fact \ref{bddfact} in general. (Hint: use Shilov boundary from \S\ref{redsec} to show that for a $k$-affinoid algebra $\calA$ with an element $f$ the spectral seminorm $\rho(f)$ can be described as the minimal number $r$ such that for any $a\in k^a$ with $|a|>r$ the element $f+a\in\calA\wtimes_k k^a$ possesses a root of any natural degree prime to $\cha(\tilk)$.)
\end{exer}

Although this fact is convenient for some applications (especially in rigid geometry), it seems to be rather accidental. We will not use it; anyway, it does not hold for trivially valued ground fields. Here is one more example when trivially valued fields require an additional care; it shows that the class of finite and admissible homomorphisms of affinoid algebras is the right analog of the class of finite homomorphisms of affine algebras.

\begin{exer}\label{admex}
(i) Show that a finite homomorphism is admissible whenever $k$ is not trivially valued, but not in general. Also, give an example of a finite bounded homomorphism which becomes non-finite after a ground field extension. (Hint: $k[T]\to k[T]$ with different norms does the job.)

(ii)* Show that a homomorphism of $k$-affinoid algebras $\phi{\colon}\calA\to\calB$ is finite admissible if and only if its ground field extension $\phi\wtimes_k K$ is finite and admissible. (Hint: the difficult case is to show the descent. First find $K_\ur$ such that the algebras become strict over $K_\ur$ and lift to the completion $L$ of $\Frac(K\otimes_k K_\ur)$. The descent from $\phi\wtimes_k L$ to $\phi\wtimes_kK_\ur$ is easy because everything is strictly affinoid, and the decent from $\phi\wtimes_kK_\ur$ to $\phi$ can be done by hands, since $K_\ur$ has a nice explicit description.)
\end{exer}

\subsection{Affinoid domains}

\subsubsection{Affinoid spectra}
In the sequel we will develop a theory of $k_H$-analytic spaces built from spectra of $k_H$-affinoid spaces (see \cite{ct}). The two extreme choices of $H$ will correspond to the general $k$-analytic spaces and strictly $k$-analytic spaces from \cite{berihes} and \cite{bernotes}. So, from now on an intermediate group $|k^\times|\subseteq H\subseteq\bfR^\times_+$ is fixed, and if not said otherwise all $k$-analytic and $k$-affinoid spaces are assumed to be $H$-strict.

\begin{rem}
When the valuation on $k$ is not trivial, it is important to develop the theory of strictly analytic spaces because it has connections to other approaches to non-archimedean geometry: formal geometry over $\kcirc$, rigid geometry and adic geometry. We prefer to develop the general $H$-strict theory because it includes both the theory of general $k$-analytic spaces and the theory of strictly $k$-analytic spaces, and hence we do not have to distinguish these two cases in some formulations. In addition, this seems to be a slightly more "accurate" approach that keeps track of the used parameters.
\end{rem}

The category of {\em $k_H$-affinoid spectra} is the category opposite to the category of $k_H$-affinoid algebras. Its objects are topological spaces of the form $\calM(\calA)$ with a $k_H$-affinoid algebra $\calA$ and morphisms are maps of the form $\calM(f)$ for bounded homomorphisms $f{\colon}\calA\to\calB$. We stress that $\calA$ is a part of the data forming the affinoid spectrum. The $k_H$-affinoid spectra are objects of global nature that will later be enriched to more geometrical affinoid spaces. This will be done in the next two sections; we will localize the current construction by introducing an appropriate Grothendieck topology and structure sheaf.

\subsubsection{Generalized normed localization}
Topology on affine schemes is defined by localization. For a $k_H$-affinoid algebra $\calA$ and an element $f\in\calA$, the localization $\calA_f$ is not affinoid for an obvious reason --- we did not worry to extend the norm. The formula $\calA_f=\calA[T]/(Tf-1)$ leads to an idea to consider the $k$-affinoid algebras $\calA_{r^{-1}f}=\calA\{rT\}/(Tf-1)$. It turns out that the latter normed localization is not general enough but its natural extension described below does the job. In the sequel, let $X=\calM(\calA)$ be a $k_H$-affinoid spectrum.

\begin{defex}\label{ratdef}
(i) Assume that elements $g,f_1\. f_n\in\calA$ do not have common zeros and $r_1\. r_n\in\sqrt{H}$ are positive numbers. Show that $$\calA_V=\calA\left\{\ur^{-1}\frac{f}{g}\right\}:=\calA\{r_1^{-1}T_1\. r_n^{-1}T_n\}/(gT_1-f_1\. gT_n-f_n)$$ is the universal $\calA$-affinoid algebra such that
$\rho_{\calA_V}(f_i)\le r_i\rho_{\calA_V}(g)$ for $1\le i\le n$. Deduce that the map $\phi_V{\colon}\calM(\calA_V)\to X$ is a bijection onto
$$V=X\left\{\ur^{-1}\frac{f}{g}\right\}:=\{x\in X|\ |f_i(x)|\le r_i|g(x)|,\ 1\le i\le n\}$$
Show that $\phi_V$ satisfies the following universal property: any morphism of $k$-affinoid spectra $\calM(\calB)\to X$ with image in $V$ factors through $\calM(\calA_V)$. The compact subset $V$ is called an $H$-strict {\em rational domain} in $X$ and by the universal property one can identify it with the $k_H$-affinoid spectrum $\calM(\calA_V)$.

(ii) For any choice of $g_1\. g_m,f_1\. f_n\in\calA$ and $s_1\. s_m, r_1 \. r_n\in\sqrt{H}$ introduce analogous domains $$V=X\{\ur^{-1}\uf,\us\ \ug^{-1}\}=\{x\in X|\ |f_i(x)|\le r_i,|g_j(x)|\ge s_j\}$$ with algebras $$\calA_V=\calA\{\ur^{-1}\uf,\us\ \ug^{-1}\}=\calA\{\ur^{-1}\uT,\us\ \uR\}/(T_i-f_i,g_jR_j-1)$$ and prove the analogs of all properties from (i). These $H$-strict domains are called {\em Laurent} (resp. {\em Weierstrass} in the case when $m=0$).

(iii) Show that any Laurent domain is rational and Laurent domains form a fundamental family of neighborhoods of a point whenever $H\neq 1$. Show that the latter is false for $H=1$.

(iv) Show that $V\subseteq X$ is an $H$-strict rational, Laurent or Weierstrass domain if and only if it is a domain of the same type and the $k$-affinoid algebra $\calA_V$ is $H$-strict. (Hint: use Fact \ref{basicaffact}.)

(v) For any map of $k$-affinoid spectra the preimage of a rational, Laurent or Weierstrass domain is a domain of the same type and given by the same inequalities. In particular, all three classes of domains are closed under finite intersections.

(vi) Show that the classes of rational and Weierstrass domains are transitive (e.g. if $Y$ is rational in $X$ and $Z$ is rational in $Y$ then $Z$ is rational in $X$), but Laurent domains are not transitive. Actually, this transitivity property is the main reason to extend the class of Laurent domains.
\end{defex}

Although, we had to consider a more general type of localizations than in the theory of affine algebra, the main difference with the theory of varieties is that affinoid domains are compact and hence have to be closed. This fact will have serious consequences when we will develop the theory of coherent sheaves.

\begin{exex}
(i) Let $X=\calM(k\{\ur^{-1}\uT\})$ be a polydisc with center at $0$ and of polyradius $\ur$, let $s_i\le r_i$ be positive numbers, and let $a_i\in k$ be elements with $|a_i|\le r_i$. Then the polydisc $\calM(k\{s_1^{-1}(T_1-a_1)\. s_n^{-1}(T_n-a_n)\})$ with center at $\ua$ and of polyradius $\us$ is a Weierstrass domain in $X$.

(ii) For $s\le r$ the annulus $A(0;s,r)=\calM(k\{r^{-1}T,sT^{-1}\})$ is a Laurent but not Weierstrass domain in the disc $E(0,r)=\calM(k\{r^{-1}T\})$.

(iii) Any finite union of discs in $E(0,r)$ is a Weierstrass domain (and is a disjoint union of finitely many discs). In particular, even when $\calA$ is an integral domain, its generalized localization does not have to be integral.
\end{exex}

\subsubsection{General affinoid domains}
It is difficult to describe a general open affine subscheme explicitly but one can easily characterize it by a universal property. Here is an affinoid analog, which was already checked for rational domains in \ref{ratdef}.

\begin{defin}
A closed subset $V\subseteq X$ is called a {\em $k_H$-affinoid domain} if there exists a morphism of $k_H$-affinoid spectra $\phi{\colon}\calM(\calA_V)\to X$ whose image coincides with $V$ and such that any morphism of $k_H$-affinoid spectra $\calM(\calB)\to X$ with image in $V$ factors through $\calM(\calA_V)$.
\end{defin}

Note that $\calA_V$ is unique up to a canonical isomorphism. The following fact allows us to identify $V$ with the $k_H$-affinoid spectrum $\calM(\calA_V)$.

\begin{fact}\label{fibfact}
The non-empty fibers of $\phi$ are isomorphisms, i.e. $\phi$ is a bijection onto $V$ and for any $y\in\calM(\calA_V)$ we have that $\calH(\phi(y))\toisom\calH(y)$.
\end{fact}

\begin{exer}\label{ggex}
(i) Prove Fact \ref{fibfact} for a point $y$ with $[\calH(y):k]<\infty$.

(ii)* Prove Fact \ref{fibfact} in general. (Hint: first prove that $\calA_V\wtimes_\calA\calA_V\toisom\calA_V$ (in particular, the separated completion homomorphism $\calA_V\otimes_\calA\calA_V\to\calA_V\wtimes_\calA\calA_V$ usually has a huge kernel); also, use without proof a non-trivial result of Gruson that the completion homomorphism $\calB\otimes_k\calB\to\calB\wtimes_k\calB$ is injective for any $k$-Banach algebra $\calB$.)

(iii)* Show that $V$ is Weierstrass if and only if the image of the homomorphism $\calA\to\calA_V$ is dense. (Hint: first you should establish the following very useful fact about affinoid generators. Assume that $\phi{\colon}k\{r_1^{-1}T_1\. r_n^{-1}T_n\}\to\calA$ is an admissible surjection, $\|\ \|$ is the residue norm on $\calA$ and $f_i=\phi(T_i)$. Then there exists $\veps>0$ such that for any choice of $f'_i\. f'_n$ with $\|f_i-f'_i\|<\veps$ there exists an admissible surjection $\phi'{\colon}k\{r_1^{-1}T_1\. r_n^{-1}T_n\}\to\calA$ with $\phi'(T_i)=f'_i$.)
\end{exer}

The following result shows that our definition of generalized localization was general enough. It allows to use rational domains for all local computations on affinoid spectra.

\begin{fact}[Gerritzen-Grauert theorem]
Any $k_H$-affinoid domain $V\subseteq X$ is a finite union of $H$-strict rational domains $V_i\subseteq X$.
\end{fact}

\begin{rem}
(i) This result (in the strict case) was not available in the first version of rigid geometry due to Tate. For this reason, Tate simply worked with rational domains and did not consider the general affinoid domains. In rigid geometry, the theorem was proved by Gerritzen-Grauert, and in Berkovich geometry the non-strict case was first deduced by Ducros. Two known rigid-theoretic proofs of this result are rather long and difficult. Originally, this theorem was needed to develop the very basics of analytic geometry, including Fact \ref{fibfact}. Later it was shown in \cite{temgg} that Fact \ref{fibfact} can be proved independently (via the hint from Exercise \ref{ggex}(ii)), and then Gerritzen-Grauert theorem can be deduced rather shortly. In addition, one obtains a similar description of all monomorphisms in the category of affinoid spaces.

(ii) To be slightly more precise, the above form of the theorem is missing in the literature. It is proved in \cite{temgg} that one can choose $V_i$'s to be rational domain. The fact that $V_i$'s may be also chosen to be $H$-strict requires a simple additional argument. For $H\neq 1$ this is done similarly to the proof of \cite[7.3]{ct}, and for $H=1$ one can deduce this from an analogous result for schemes (because for a trivially valued $k$, the category of $k_1$-affinoid algebras is equivalent to the category of finitely generated $k$-algebras.)
\end{rem}

\subsection{$G$-topology and the structure sheaf}\label{Gsec}
From now on we assume that $H\neq 1$.

\subsubsection{$G$-topology of compact domains}
In order to define $k$-affinoid spaces we should provide each spectrum $X=\calM(\calA)$ with a certain structure sheaf $\calO_X$. Naturally, we would like $\calO_X$ to be a sheaf of $k$-affinoid or $k$-Banach algebras but then we should study the sections over closed subsets (e.g. affinoid domains). A naive attempt to consider the topology generated by affinoid domains does not work out.

\begin{exer}
Observe that the unit interval $[0,1]$ is neither connected nor compact in the topology generated by closed intervals $[a,b]$. Show, similarly, that the unit closed disc $\calM(k\{T\})$ is neither connected nor compact in the topology generated by affinoid domains.
\end{exer}

A brilliant idea of Tate (with a strong influence of Grothendieck) is to generalize the notion of topology by allowing only certain open coverings. The resulting notion of a $G$-topology $\tau$ is simply a Grothendieck topology on a set $\tau^\opp$ of subsets of $X$ such that $\tau^\opp$ is closed under finite intersections and any covering of this topology is also a set-theoretical covering. Sets of $\tau^\opp$ are also called $\tau$-open or and the coverings of this Grothendieck topology are called {\em admissible coverings}. (Note that $X$ is in $\tau^\opp$ because it is the intersection indexed by the empty set.)

\begin{defin}
(i) A {\em compact $k_H$-analytic domain} $Y$ in a $k_H$-affinoid spectrum $X$ is a finite union of $k_H$-affinoid domains. (It is called a special domain in \cite{berbook}.)

(ii) The {\em $H$-strict compact $G$-topology} $\tau^c_H$ on a $k_H$-affinoid spectrum $X$ is defined as follows: $(\tau^c_H)^\opp$ is the set of all compact $k_H$-analytic domains and admissible coverings are the finite ones.
\end{defin}

\begin{rem}
By Gerritzen-Grauert theorem one can replace affinoid domains with rational domains in this definition obtaining the original Tate's definition of $G$-topology.
\end{rem}

\subsubsection{The structure sheaf}
Tate proved that (in the strict case) this $G$-topology is the right tool to define coherent sheaves of modules. In particular, $\calO_{X_H}(V)=\calA_V$ extends to a $\tau^c_H$-sheaf of Banach algebras.

\begin{fact}[Tate's acyclity theorem]\label{tatefact}
For any finite affinoid covering $X=\cup_i V_i$ and finite Banach $\calA$-module $M$ the \v{C}hech complex $$0\to M\to\prod_i M_i\to\prod_{i,j} M_{ij}\to\dots$$
is exact and admissible, where $M_i=M\otimes_\calA\calA_{V_i}$, $M_{ij}=M\otimes_\calA\calA_{V_{ij}}$, etc.
\end{fact}

Admissibility in this result is very important. In particular, it allows to define norms on the structure sheaf introduced below.

\begin{exer}
(i) For any compact $k_H$-analytic domain $V$ with a finite affinoid covering $V=\cup_i V_i$ set $\calO_{X_H}(V)=\Ker(\prod_i \calA_{V_i}\to\prod_{i,j}\calA_{V_i\cap V_j})$ and provide it with the restriction norm. Show that $\calO_{X_H}(V)$ depends only on $V$, in particular, $\calO_{X_H}(V)=\calA_V$ when $V$ is affinoid. In addition, show that $\calO_{X_H}$ is a sheaf of $k$-Banach algebras, in particular, the restriction morphisms are bounded.

(ii) Deduce that any polydisc $X=\calM(k\{r_1^{-1}T_1\.r_n^{-1}T_n\})$ with $r_i\in\sqrt{H}$ is
$\tau^c_H$-connected, i.e. $X$ is not a disjoint union of two non-empty compact $k_H$-analytic domains. (Hint: $\calO_{X_H}(X)$ is integral.)
\end{exer}

\begin{fact}\label{afffact}
A compact $k_H$-analytic domain $V\subseteq X=\calM(\calA)$ is an affinoid domain if and only if the Banach algebra $\calA_V=\calO_{X_H}(V)$ is $k_H$-affinoid and the natural map of sets $\phi_V{\colon}V\to\calM(\calA_V)$ is bijective. In particular, any $k_H$-affinoid domain $V$ is a $k$-affinoid domain.
\end{fact}

\begin{exer}
(i) Define the map $\phi_V$ in Fact \ref{afffact}. (Hint: take a finite affinoid covering $V=\cup_{i=1}^n V_i$ and show that the affinoid morphisms $V_i\to\calM(\calA_V)$ are compatible on intersections.)

(ii)* Prove Fact \ref{afffact}. (Hint: use Gerritzen-Grauert and Tate's acyclity theorems; the main stage is to show that if $U=X\{\ur^{-1}\uf/g\}$ is a rational domain contained in $V$ then $U\toisom\calM(\calA_V\{\ur^{-1}\uf/g\})$.)
\end{exer}

\subsubsection{$k_H$-affinoid spaces}
\begin{defin}
(i) A {\em $k_H$-affinoid space} $X$ is a $k_H$-affinoid spectrum $\calM(\calA)$ provided with the $G$-topology $\tau^c_H$ and the $\tau^c_H$-sheaf of $k$-Banach rings $\calO_{X_H}$, which is called the {\em structure sheaf}.

(ii) A {\em morphism} of $k_H$-affinoid spaces is a continuous and $\tau^c_H$-continuous map $f{\colon}Y\to X$ provided with a {\em bounded} homomorphism of sheaves $f^\#{\colon}\calO_{X_H}\to f_*(\calO_{Y_H})$ in the sense that for any pair of compact $k_H$-analytic domains $X'\subseteq X$ and $Y'\subseteq f^{-1}(X)$ the homomorphism
$f^\#{\colon}\calO_{X_H}(X')\to\calO_{Y_H}(Y')$ is bounded.

(iii) A $\tau^c_H$-sheaf of finite $\calO_{X_H}$-Banach modules is {\em coherent} if it is of the form $\calM(V)=M\otimes_\calA\calO_{X_H}(V)$ for a finite Banach $\calA$-module $M$.
\end{defin}

\begin{exer}
Show that the continuity assumption in (ii) follows from $\tau^c_H$-continuity and hence can be removed.
\end{exer}

Note that $k_H$-affinoid spaces are "self-contained" geometric objects analogous to affine schemes. It will later be an easy task to globalize this definition.

\begin{fact}\label{newfact}
The categories of $k_H$-affinoid spectra and $k_H$-affinoid spaces are naturally equivalent.
\end{fact}

Since this fact seems to be new, we give a detailed exercise on its proof.

\begin{exer}
(i) Reduce the fact to the following claim: if $(f,f^\#){\colon}(Y,\calO_{Y_H})\to(X,\calO_{X_H})$ is a
morphism of $k_H$-affinoid spaces and $\phi{\colon}\calO_{X_H}(X)\to\calO_{Y_H}(Y)$ is the induced bounded homomorphism of $k_H$-affinoid algebras, then $f=\calM(\phi)$ and $f^\#{\colon}\calO_{X_H}\to f_*(\calO_{Y_H})$ is the bounded homomorphism of the structure sheaves induced by $\phi$. In other words, $(f,f^\#)$ is, in its turn, induced by $\phi$.

(ii) Choose a point $y\in Y$ with $x=f(y)$. For any $H$-strict rational domain $X'=X\{r^{-1}\frac{g}{h}\}$ containing $x$ choose an $H$-strict rational domain $Y'\subseteq f^{-1}(X')$ containing $y$. Set $\calA=\calO_{X_H}(X)$, $\calB=\calO_{Y_H}(Y)$ and $\calB'=\calO_{Y_H}(Y')$. Since the homomorphism $\calA\to\calB'$ factors through $\calO_{X_H}(X')=\calA\{r^{-1}\frac{g}{h}\}$, it follows that the homomorphism
$\calB\to\calB'$ factors through $\calB\wtimes_\calA\calA\{r^{-1}\frac{g}{h}\}\toisom\calB\{r^{-1}\frac{g}{h}\}$, and hence $Y'\subseteq Y\{r^{-1}\frac{g}{h}\}$. Therefore, any inequality $|g(x)|\le r|h(x)|$ with $g,h\in\calO_{X_H}(X)$ implies that $|\phi(g)(y)|\le r|\phi(h)(y)|$. Deduce that $\calM(\phi)$ takes $y$ to $x$ and hence coincides with $f$.

(iii) Finish the argument by showing that $f^\#$ is also induced by $\phi$. (Hint: check this for sections on rational domains and then apply Tate's theorem.)
\end{exer}

In the sequel we will not distinguish between $k_H$-affinoid spectra and $k_H$-affinoid spaces, that is, we will automatically enrich any $k_H$-affinoid spectrum with the structure of a $k_H$-affinoid space. Also, we will refine the structure of $k_H$-affinoid spaces a little bit more in \S\ref{affansec}.

\subsubsection{Coherent sheaves}
We finish \S\ref{Gsec} with a discussion on coherent sheaves. By Tate's theorem any finite Banach $\calA$-module gives rise to a sheaf of finite Banach $\calO_{X_H}$-modules. The opposite result (in a slightly different formulation) was proved by Kiehl.

\begin{fact}
(i) (Kiehl's theorem) Any $G$-locally coherent sheaf is coherent. Namely, if for a finite Banach $\calO_{X_H}$-module $\calM$ there exists a finite affinoid covering $X=\cup_i V_i$ such that the restrictions $\calM|_{V_i}$ are coherent then $\calM$ is coherent.

(ii) Tate's and Kiehl's theorems easily imply that the categories of coherent $\calO_{X_H}$-module and finite Banach $\calA$-modules are naturally equivalent.
\end{fact}

\begin{rem}
The theory of $k_H$-analytic (and even $k_H$-affinoid) spaces does not admit a reasonable class of infinite type modules analogous to the class of quasi-coherent modules on schemes. Moreover, this theory does not even have a notion of affinoid morphisms. There exist morphisms $f{\colon}Y\to X$ with finite affinoid coverings $X=\cup_i X_i$ such that $X$ and $f^{-1}(X_i)$ are affinoid but $Y$ is not affinoid, see Example \ref{nogoodexam}(ii).
\end{rem}

\subsection{The reduction map, boundary and interior}

\subsubsection{Reduction}\label{redsec}
In general, reduction relates (strictly) $k$-affinoid algebras and spaces to geometry over the residue field $\tilk$.

\begin{fact}
Assume that the valuation is non-trivial. For a bounded homomorphism $\phi{\colon}\calA\to\calB$ of strictly $k$-affinoid algebras the following conditions are equivalent: $\phi$ is finite,
$\tilphi{\colon}\tilcalA\to\tilcalB$ is finite, $\phicirc{\colon}\calAcirc\to\calBcirc$ is integral.
\end{fact}

\begin{exer}
Deduce that for any strictly $k$-affinoid $\calA$ the reduction $\tilcalA$ is finitely generated over $\tilk$.
\end{exer}

\begin{rem}
(i) This result shows that the reduction functor controls strictly $k$-affinoid algebras very well. A similar result holds for general $k_H$-affinoid algebras if one replaces $\tilA$ with the $H$-graded reduction $$\tilA_H=\oplus_{h\in H}\{x\in\calA|\ \rho(x)\le h\}/\{x\in\calA|\ \rho(x)<h\}$$

(ii) The question whether $\phicirc$ is finite is more subtle. Already for a finite field extension $l/k$ one often has that $\lcirc/\kcirc$ is not finite. On the other hand if $K$ is algebraically closed or discretely valued and the algebras are reduced then $\phicirc$ is integral if and only if it is finite. We refer the reader to \cite[Ch 6, \S\S3-4]{BGR} for more details.
\end{rem}

Now, let us study the geometric side of the reduction.

\begin{defin}
(i) {\em Reduction} of a strictly $k$-affinoid space $X=\calM(\calA)$ is the reduced $\tilk$-variety $\tilX=\Spec(\tilcalA)$.

(ii) {\em Reduction map} $\pi_X{\colon}X\to\tilX$ sends a point $x$ with the character
$\chi_x{\colon}\calA\to\calH(x)$ to the point $\tilx\in\tilX$ induced by the character
$\tilchi_x{\colon}\tilcalA\to\wHx$. (Note that $k(\tilx)$ can be much smaller than $\wHx$; as we know from Exercise \ref{typesex}(5) the latter field does not even have to be finitely generated over $\tilk$.)
\end{defin}

\begin{exer}
The map $\pi_X$ is {\em anti-continuous} in the sense that the preimage of an open set is closed and vice versa.
\end{exer}

\begin{fact}
(i) The reduction map is surjective.

(ii) The preimage of a generic point of $\tilX$ is a single point, and the union of all such points is the {\em Shilov boundary} $\Gamma(X)$ of $X$. Namely, any function $|f(x)|$ with $f\in\calA$ accepts its maximum on $\Gamma(X)$ and $\Gamma(X)$ is the minimal closed set satisfying this property.
\end{fact}

\begin{rem}\label{grrem}
The same result holds for $H$-graded reduction if one defines $\tilX_H=\Spec_H(\tilcalA_H)$ as the set of all homogeneous prime ideals in the $H$-graded ring $\tilcalA_H$.
\end{rem}

\begin{exam}
(i) The spectral seminorm on $\calA$ is multiplicative if and only if $\tilcalA$ is integral. In this case, the spectral seminorm itself defines a point which is both the preimage of the generic point of $\tilX$ and the Shilov boundary of $X$. For example, the Shilov boundary of a polydisc $E(0,\ur)=\calM(k\{\ur^{-1}\uT\})$ is a single (maximal) point.

(ii) The Shilov boundary of a closed annulus $X=A(0;s,r)=E(0,r)\setminus D(0,s)=\calM(k\{r^{-1}T,sT^{-1}\})$ with $s<r$ consists of two points, the maximal points of the closed discs $E(0,r)$ and $E(0,s)$. Assume that $s,r\in |k^\times|$, then the reduction $\tilX$ is the cross $\Spec(\tilk[R,S]/(RS))$, where $R$ and $S$ are the reductions of appropriate rescalings of $T$ and $T^{-1}$. The points collide when $s$ tends to $r$, namely, $\Gamma(A(0;r,r))$ consists of a single point and the reduction is $\tilk[R,S]/(RS-1)=\tilk[\tilT,\tilT^{-1}]$ in this case.
\end{exam}

\subsubsection{Relative boundary and interior}\label{boundsec}

\begin{defin}
Let $\phi{\colon}Y\to X$ be a morphism of $k$-affinoid spaces and let $X=\calM(\calA)$ and $Y=\calM(\calB)$. {\em Relative interior} $\Int(Y/X)\subseteq Y$ consists of points $y\in Y$ such that there exists an admissible surjective $\calA$-homomorphism $\psi{\colon}\calA\{r_1^{-1}T_1\.r_n^{-1}T_n\}\to\calB$ with $|\psi(T_i)(y)|<r_i$ for $1\le i\le n$. {\em Relative boundary} $\partial(Y/X)$ is the complement to the relative interior in $Y$.  Absolute {\em interior} $\Int(Y)$ and {\em boundary} $\partial(Y)$ are defined with respect to the morphism $Y\to\calM(k)$. A morphism
\end{defin}

\begin{rem}
(i) This definition has a very geometric interpretation as follows. The homomorphism $\psi$ induces a closed immersion of $Y$ into a relative closed polydisc of polyradius $\ur$ over $X$. The inequalities in the definition mean that the image of $y$ lies in the open relative polydisc of the same polyradius.

(ii) Realization of boundary as a set is specific to Berkovich analytic geometry. For example, we will later see that boundaryless morphisms of general spaces form a very important class (e.g. any proper morphism is boundaryless). This class was introduced in the classical rigid geometry in a rather formal way, but it was not so clear how to describe obstructions for being boundaryless (note that analytic boundary consists of non-rigid points). To some extent this obstruction could be felt by working with formal models, but a concrete notion of boundary was missing.
\end{rem}

The notions of relative boundary and interior turn out to be very important in analytic geometry. Most of the facts about them are proved by use of the reduction theory. Here is a list of their basic properties.

\begin{fact}\label{intboufact}
(i) The relative interior is open in $Y$ and the relative boundary is closed.

(ii) The relative interiors are compatible with compositions in the sense that
$\Int(Z/X)=\Int(Z/Y)\cap\psi^{-1}(\Int(Y/X))$ for a pair of morphism $Z\stackrel{\psi}{\to}Y\stackrel{\phi}{\to}X$.

(iii) Relative boundary is $G$-local on the base in the sense that for a finite affinoid covering $X=\cup_i X_i$ and $Y_i=\phi^{-1}(X_i)$ one has that $\partial(Y/X)=\cup_i\partial(Y_i/X_i)$.

(iv) $\phi$ is {\em boundaryless}, i.e. has an empty boundary, if and only if it is {\em finite}, i.e. $\calA\to\calB$ is finite admissible.

(v) If $Y$ is an affinoid domain in $X$ then $\Int(Y/X)$ is the topological interior of $Y$ in $X$.

(vi) Assume that $X$ and $Y$ are strictly $k$-affinoid, $y\in Y$ and $x=\phi(y)$. Then $y\in\Int(Y/X)$ if and only if the image of $\tilcalB$ in $k(\tily)$ (usually denoted $\tilchi_y(\tilcalB)$) is finite over the image of $\tilcalA$ in $k(\tilx)$.
\end{fact}

\begin{rem}
The last condition is very convenient for explicit computations. Its analog holds for $H$-strict $X$ and $Y$ and $H$-graded reduction.
\end{rem}

Now, we illustrate the introduced notions with some examples. For simplicity, all analytic spaces in the examples are assumed to be strict.

\begin{exex}
(i) Show that $\Int(Y)$ is the preimage under the reduction map $\pi_Y$ of the set of closed points of $\tilY$.

(ii) Let $X$ be an {\em affinoid curve}, i.e. a $k$-affinoid space of dimension one (in the sense of \S\ref{dimsec} below) and without isolated Zariski closed points. Show that the boundary of $X$ coincides with its Shilov boundary. (Hint: use Noether normalization to prove that $\tilX$ is a curve.)
\end{exex}

Next, let us consider a higher dimensional example.

\begin{exex}
Let $\phi$ be the projection of the unit polydisc $Y=\calM(k\{T,S\})$ onto the unit disc $X=\calM(k\{S\})$. We will compare $\partial(Y)$ and $\partial(Y/X)$ fiberwise over $X$. Let $x\in X$ be a point and let $Y_x=\psi^{-1}(x)$ denote the fiber over $x$, so that $Y_x=\calM(\calH(x)\{T\})=E_{\calH(x)}(0,1)$ is a closed unit disc over $\calH(x)$.

(i) The fiber disc $Y_x$ contains exactly one point in its boundary $\partial(Y_x)$, the maximal point of the disc.

(ii) If $x$ is not the maximal point of $X$ then the sets $Z=\partial(Y/X)\cap Y_x$ and $\partial(Y)\cap Y_x$ coincide by Fact \ref{intboufact}(ii). A non-maximal point of $Y_x$ is contained in $Z$ if and only if it is contained in a subdisc $E_{\calH(x)}(a,r)$ of $Y_x$
with $r<1$ and $\inf|a-k^a|=1$. In particular, if $x$ is of type 1,3 or 4 then there are no such points and so $Z=\partial(Y_x)$, but if $x$ is of type 2 then $Z$ contains infinitely many open unit subdiscs.

(iii) If $x$ is the maximal point of $X$ then $\partial(Y/X)\cap Y_x$ is much smaller than $\partial(Y)\cap Y_x$. For example, let $\eta$ and $\ve$ be the generic points of the quadrics $\tilT^2-\tilS=0$ and $\tilT\tilS=1$ in $\tilY=\Spec(\tilk[\tilT,\tilS])$. Show that $\pi_Y^{-1}(\eta)\subset\Int(Y/X)$ but $\pi_Y^{-1}(\ve)\subset\partial(Y/X)$. Find a geometric explanation for this fact. (Hint: how are these fibers embedded into a larger polydisc $\calM(k\{r^{-1}T,S\})$ with $r>1$?)
\end{exex}

\subsection{The dimension theory}\label{dimsec}
\subsubsection{} In the strict case one can just copy the dimension theory of affine varieties. However, we will see that in general one has to be slightly more careful.

\begin{exer}
(i) Prove that the dimension of a strictly $k$-affinoid algebra $\calA$ is preserved under any ground field extension. (Hint: use Noether normalization.)

(ii) Show that general $k$-affinoid algebras do not share this property. (Hint: show that $K_r\wtimes_k K_r\toisom K_r\{r^{-1}T\}\toisom K_r\{T\}$.)
\end{exer}

Since the dimension stabilizes after a sufficiently large ground field extension, it is natural to define the dimension of $k$-affinoid spaces as follows.

\begin{defin}
Dimension $\dim(X)$ of a $k$-affinoid space $X=\calM(\calA)$ is the dimension of an algebra
$\calA_K=\calA\wtimes_k K$, where $K$ is a non-archimedean $k$-field such that $\calA_K$ is strictly $K$-affinoid.
\end{defin}

\begin{rem}
Assume for concreteness that $r=(r_1)$. One can view $X=\calM(K_r)$ as a $k$-curve consisting of its generic point. The only difference with the theory of algebraic $k$-curves is that $X$ is of "finite type" over $k$.
\end{rem}

The following exercise illustrates that type 2, 3 and 4 points are sort of "generic points" of the curves, so (informally) one can imagine them as points of dimension 1.

\begin{exer}
Let $x$ be a point of $\bfA^1_k$.

(i) Show that $x$ is of type 1 if and only if the fiber over $x$ in any ground field extension $\bfA^1_K$ is profinite.

(ii) Show that $x$ is not of type 1 if and only if the fiber over $x$ in some $\bfA^1_K$ contains a closed $K$-disc.
\end{exer}

\section{Analytic spaces}\label{4}

\subsection{The category of $k$-analytic spaces}\label{affansec}

\subsubsection{Nets}

\begin{defin}
Let $X$ be a topological space with a set of subsets $T$.

(i) $T$ is a {\em quasi-net} if any point $x\in X$ has a neighborhood of the form $\cup_{i=1}^n V_i$ with $x\in V_i\in T$ for $1\le i\le n$.

(ii) A quasi-net $T$ is a {\em net} if for any choice of $U,V\in T$ the restriction $T|_{U\cap V}=\{W\in T|\ W\subseteq U\cap V\}$ is a quasi-net of subsets of the topological space $U\cap V$.
\end{defin}

We remark that the definition of nets axiomatizes properties an affinoid atlas should satisfy, and the definition of quasi-nets axiomatizes the properties of admissible coverings by analytic domains (see Remark \ref{admrem} below). Now, let us assume that $X$ is locally Hausdorff and the sets of $T$ are compact (as will be the case with analytic spaces and their atlases).

\begin{exer}\label{netex}
(i) With the above assumptions, $X$ is locally compact.

(ii) A subset $Y\subseteq X$ is open if and only if $Y\cap U$ is open in $U$ for any $U\in T$.

(iii) $X$ is Hausdorff if and only if for any pair $U,V\in T$ the intersection $U\cap V$ is compact.

(iv) A subset $Y\subseteq X$ is compact if and only if it is covered by compact intersections of the form $Y\cap U$ with $U\in T$ (but there may exist non-compact intersections with other elements of $T$).
\end{exer}

\begin{rem}
Nets in our sense are slightly analogous to $\ve$-nets on metric spaces, and they should not be confused with the nets that generalize sequences in the definition of limits.
\end{rem}

\subsubsection{Analytic spaces}
We will freely view a net as a category with morphisms being the inclusions.

\begin{defin}
A {\em $k_H$-analytic space} is a locally Hausdorff topological space $X$ with an {\em atlas of $k_H$-affinoid domains}, where the atlas consists of a net $\tau_0$ on $X$, a functor $\phi{\colon}\tau_0\to k_H$-$Aff$ taking morphisms in $\tau_0$ to embeddings of affinoid domains, and an isomorphism $i$ of the two natural topological realization functors from $\tau_0$ to the category of topological spaces, $Top$ and $Top\circ\phi$. In concrete terms, we will write $\phi(U)=\calM(\calA_U)$, and $i$ reduces to giving a homeomorphism $i_U{\colon}U\toisom\phi(U)$ for any $U\in\tau_0$ such that for any inclusion $j{\colon}U\into V$ in $\tau$ these homeomorphisms are compatible with $\phi(j){\colon}\phi(U)\to\phi(V)$.
\end{defin}

By Exercise \ref{netex}(i) $X$ is locally compact.

\begin{defin}
(i) A {\em $k_H$-analytic domain} in $X$ is a subset $V\subseteq X$ that admits a covering $V=\cup_{i\in I}V_i$ such that each $V_i$ is a $k_H$-affinoid domain in some element of $\tau_0$ and $\{V_i\}_{i\in I}$ is a quasi-net on $V$.

(ii) By $\tau_H$ (resp. $\tau^c_H$) we denote the sets of all (resp. compact) $k_H$-analytic domains. A covering of an element $V\in\tau_H$ by elements $V_i\in\tau_H$ is {\em admissible} if $V_i$'s form a quasi-net on $V$.
\end{defin}

\begin{rem}\label{admrem}
Note that a closed unit disc $X$ is a disjoint union of the open unit disc $D(0,1)$ and the closed annulus $A(0;1,1)$. However, the covering $X=D(0,1)\coprod A(0;1,1)$ is not admissible thanks to the condition $x\in\cap_{i=1}^nV_i$ in the definition of quasi-nets. In particular, one can easily show that $X$ is $\tau_H$-connected. This explains the role of the first condition in the definition of quasi-nets, and the second condition is needed to ensure that analytic spaces are locally compact.
\end{rem}

\begin{exer}
(i) Show that $\tau_H$ with admissible coverings is a $G$-topology and give an example where it is not closed under finite unions. In the sequel we will refer to $\tau_H$ as the {\em $G$-topology of $X$}. (Hint: already the union of the open polydisc of polyradius $(1,2)$ with the closed polydisc of polyradius $(2,1)$ is not locally compact.)

(ii) Show that although $\tau^c_H$ does not have products in general (e.g. $X$ is not in $\tau^c_H$ if it is not compact and $\tau^c_H$ is not closed under intersections of pairs for a non-Hausdorff space), it is closed under fibred products (i.e. if $U,V\subset W$ are three elements of $\tau^c_H$ then $U\cap V$ is in $\tau^c_H$). Use this to define $\tau^c_H$-sheaves.

(iii) Show that any admissible covering of a compact $k_H$-analytic domain by compact $k_H$-analytic domains possesses a finite subcovering. Deduce that $V\in\tau^c_H$ if and only if it $V=\cup_{i=1}^nV_i$ with each $V_i$ being $k_H$-affinoid in an element of $\tau_0$ and all intersections $V_i\cap V_j$ being compact.

(iv) Deduce that the correspondence $V\to\calA_V$ on $\tau_0$ extends uniquely to a $\tau^c_H$-sheaf of Banach $k$-algebras. (Hint: find an affinoid covering as in (iii), set $\calA_V=\Ker(\prod_i\calA_{V_i}\to\prod_{i,j}\calA_{V_{ij}})$, and use Tate's acyclity theorem to establish independence of the covering.)
\end{exer}

The sheaf from (iv) will be called the {\em structure sheaf} and denoted $\calO_{X_H}$. Any $k_H$-analytic space given by an atlas will automatically be provided with this additional structure.

\subsubsection{Morphisms between analytic spaces}
Intuitively, a morphism between $k_H$-analytic spaces should be a continuous and $G$-continuous map $f{\colon}Y\to X$ (i.e. the preimage of an analytic domain is an analytic domain) provided with a bounded homomorphism $f^\#{\colon}\calO_{X_H}\to f_*\calO_{Y_H}$. However, the
direct image $f_*\calO_{Y_H}$ does not really make sense for non-compact morphisms, including
$\bfA^1_k\to\calM(k)$. Therefore, we suggest the following definition.

\begin{defin}
A morphism $f{\colon}Y\to X$ between $k_H$-analytic spaces consists of a continuous and $G$-continuous map $f{\colon}Y\to X$ and a compatible family of bounded homomorphisms
$f^\#_{U,V}{\colon}\calO_{X_H}(U)\to\calO_{Y_H}(V)$ for any pair of compact $k_H$-analytic domains $U\subseteq X$ and $V\subseteq f^{-1}(U)$.
\end{defin}

A priori it is not clear how to compose such morphisms because the image of a compact set does not have to be Hausdorff. This forces us to show that a morphism is determined already by its restriction to atlases. (Note that the atlas definition of morphisms is used in \cite[\S1.2]{berihes} and \cite[3.1]{bernotes}.)

\begin{exer}
(i) Assume that $Y$ and $X$ are provided with affinoid atlases $\tau_Y$ and $\tau_X$ such that for any $U\in\tau_X$ the restriction $\tau_Y|_{f^{-1}(U)}$ is an atlas of $f^{-1}(U)$. Then to give a morphism $Y\to X$ is equivalent to give a similar data $(g,g^\#)$but with $g^\#_{U,V}$ defined only for $U\in\tau_X$ and $V\in\tau_Y$ with $f(V)\subseteq U$.

(ii) Use this to define composition of morphisms.

(iii) Show that any $k_H$-analytic domain $V\subseteq X$ possesses a canonical structure of a $k_H$-analytic space. Moreover, the inclusion underlies the {\em analytic domain embedding} morphism $i_V{\colon}V\to X$ which possesses the universal property that any morphism $Y\to X$ with set-theoretical image in $V$ factors through $V$ (in the analytic category).
\end{exer}

Thanks to the claim of (ii) we have now introduced the category of $k_H$-analytic spaces, which will be denoted $k_H$-$An$. The particular case of (resp. {\em strictly}) {\em $k$-analytic spaces} corresponding to $H=\bfR^\times_+$ (resp. $H=|k^\times|$) will be denoted $k$-$An$ (resp. $st$-$k$-$An$).

\begin{rem}
(i) Note that an isomorphism between two analytic spaces $X$ and $Y$ is a homeomorphism $f{\colon}Y\to X$ which induces a bijection $\tau_{H,Y}\toisom\tau_{H,X}$ and an isomorphism of the structure sheaves $f^\#{\colon}\calO_{X_H}\to f_*(\calO_{Y_H})$. In particular, this allows to consider $k_H$-analytic spaces without any fixed affinoid atlas.

(ii) We gave the usual definition of analytic spaces that follows \cite{berihes} but used a different definition of morphisms. Berkovich first defines a category of spaces with a fixed atlas (and morphisms between such objects) and then inverts morphisms corresponding to refinements of atlases.
\end{rem}

One can also define analytic spaces in a more invariant way that does not involve atlases. This is worked out in the following exercise.

\begin{exer}
Show that the following definition of analytic spaces is equivalent to the standard one: a $k_H$-analytic space $X$ is a topological space $|X|$ provided with a $G$-topology $\tau^c_H$ and a $\tau^c_H$-sheaf of Banach $k$-algebras $\calO_{X_H}$ such that the elements of $\tau^c_H$ are compact (in the usual topology), $\tau^c_H$ is closed under finite unions, and $X$ is $G$-locally isomorphic to $k_H$-affinoid sets, in the sense that there is a quasi-net $\{X_i\}$ on $X$ such that each triple $(X_i,\tau_H|_{X_i},\calO_{X_H}|_{X_i})$ is isomorphic to a $k_H$-affinoid space. (Hint: use Fact \ref{newfact}.)
\end{exer}

The following useful result is surprisingly difficult (for non-separated spaces). It was proved in \cite{temred2} for the strictly analytic category and was generalized to general $H$-strict spaces in \cite{ct}.

\begin{fact}\label{faithfact}
Assume that $H'\subseteq\bfR_+^\times$ is a subgroup containing $H$. The natural embedding functor $k_H$-$An\to k_{H'}$-$An$ is fully faithful. In particular, any analytic space admits at most one (up to an isomorphism) structure of an $H$-strict analytic space.
\end{fact}

An equivalent way to reformulate this fact is by saying that any $k$-analytic morphisms between $H$-strict $k$-analytic spaces can be described using $H$-strict atlases. The fact implies that even when studying $k_H$-analytic spaces we can (and in the sequel will) safely work within the category of all $k$-analytic spaces. In the sequel, our default $G$-topology is the $G$-topology $\tau$ of all $k$-analytic domains, and $\tau^c$ denotes the $G$-topology of compact $k$-analytic domains.

\begin{exer}
Show that the $\tau^c$-sheaf $\calO_{X_{\bfR^\times_+}}$ uniquely extends to a $G$-sheaf of rings $\calO_{X_G}$ and, moreover, $\calO_{X_G}(V)=\Mor_{k-An}(V,\bfA^1_k)$.
\end{exer}

From now on, the {\em structure sheaf} of $X$ refers to the sheaf $\calO_{X_G}$.

\begin{rem}
(i) We aware the reader that there is an abuse of language in our notion of the structure sheaf because $\calO_{X_G}$ is not a part of the definition of $X$, but only an additional structure
$X$ is provided with. Moreover, in a sharp contrast with the $\tau^c$-sheaf $\calO_{X_{\bfR^\times_+}}$, the $G$-sheaf $\calO_{X_G}$ is not a sheaf of Banach rings and it does not contain enough information to define the analytic space (at least when the valuation is trivial).

(ii) One can provide $\calO_{X_G}$ with a natural structure of {\em pluri-normed} $k$-algebras, i.e. $k$-algebras with a family of $k$-bounded norms. If the family is countable, then this agrees with the usual notion of a Frechet $k$-algebra. Probably, one can develop a theory that includes also spectra of pluri-normed algebras, but this was not done so far (to the best of my knowledge).
\end{rem}

\subsubsection{Gluing of analytic spaces}
There are three main constructions of new $k$-analytic spaces: by gluing (or using atlases), by analytification of algebraic $k$-varieties, and as the generic fiber of formal $\kcirc$-schemes. Here we consider only the first construction because the other two will be studied later.

\begin{exer}
Assume that $\{X_i\}_{i\in I}$ is a family of $k$-analytic spaces provided with analytic domains $X_{ij}\into X_i$ and isomorphisms $\phi_{ij}{\colon}X_{ij}\toisom X_{ji}$ that satisfy the usual cocyle compatibility condition on the intersections $X_{ijk}=X_{ij}\cap X_{ik}$. Show that in the following cases they glue to a $k$-analytic space $X$ covered by domains isomorphic to $X_i$ so that $X_i\cap X_j\toisom X_{ij}$.

Case 1: The domains $X_{ij}\subseteq X_i$ are open.

Case 2; The domains $X_{ij}\subseteq X_i$ are closed and for each $i$ only finitely many domains $X_{ij}$ are non-empty.
\end{exer}

\begin{exer}\label{glueex}
(i) Let $X$ be glued from annuli $X_1=A(0;r,s)$ and $X_2=A(0;s,t)$ along $X_{12}=A(0;s,s)$ so that the orientation of the annuli is preserved. The latter means that the gluing homomorphism
$k\{s^{-1}T_1,rT_1^{-1}\}\to k\{s^{-1}T,sT^{-1}\}$ takes $T_1$ to an element $\sum_{i=-\infty}^\infty a_iT^i$ with $|a_1T|=s>|a_iT^i|$ for $i\neq 1$, and similarly for the second chart. Show that $X$ is isomorphic to the annulus $A(0;r,t)$. (Hint: show that the intersection of $k\{s^{-1}T_1,rT_1^{-1}\}$ and $k\{t^{-1}T_2,sT_2^{-1}\}$ inside of $k\{s^{-1}T,sT^{-1}\}$ is isomorphic to $k\{t^{-1}R,rR^{-1}\}$.)

(ii) In the same way show that if $X_1$ is the disc $E(0,s)$ and $X_2,X_{12}$ are as in (i)  then $X\toisom E(0,t)$.

(iii) We define $\bfP^1_k$ as the obvious gluing of $\calM(k\{T\})$ and $\calM(k\{T^{-1}\})$ along $\calM(k\{T,T^{-1}\})$. Show that any other gluing of two discs with the same choice of orientation is isomorphic to $\bfP^1_k$. A wrong choice of orientation leads to a space that we call a closed unit disc with doubled open unit disc. This space is Hausdorff, but we will later see that it is not locally separated at the maximal point of the disc.

(iv) Define $\bfP^n_k$ with homogeneous coordinates $T_0\. T_n$ in two different ways: (a) as a gluing of $n+1$ unit polydiscs, (b) as a gluing of $n+1$ affine spaces $\bfA^n_k$. (Hint: in both cases, it is convenient to symbolically denote coordinates on the $i$-th chart as $\frac{T_j}{T_i}$ for $0\le j\le n,j\neq i$.)
\end{exer}

\begin{defin}
(i) A seminorm on a graded ring $A=\oplus_{d\in\bfN}A_d$ is {\em homogeneous} if it is determined by its values on the homogeneous elements via the max formula $|\sum_{d\in\bfN} a_d|=\max_{d\in\bfN}|a_d|$.

(ii) Seminorms $|\ |$ and $\|\ \|$ on $A$ are {\em homothetic} if there exists a number $C>0$ such that $|a_d|=C^d\|a_d\|$ for any $a_d\in A_d$.

(iii) Let $k$ be a non-archimedean field with a graded $k$-algebra $A$. The projective spectrum $\MProj(A)$ is the set of all homothety equivalence classes of homogeneous semivaluations on $A$ that extend the valuation of $k$ and do not vanish on the whole $A_1$.
\end{defin}

\begin{exer}
(i) Show that $\bfP^n_k=\MProj(k[T_0\. T_n])$ by a direct computation.

(ii) Alternatively, show that two points $x,y\in\bfA^{n+1}\setminus\{0\}=\MSpec(k[\uT])\setminus\{0\}$ are mapped to the same point by the projection $\bfA^{n+1}\setminus\{0\}\to\bfP^n_k$ if and only if there exists $C>0$ such
that $|f_d(x)|=C^d|f_d(y)|$ for any homogeneous $f_d(\uT)\in k[\uT]$ of degree $d$. Show that the set of homogeneous semivaluations in any fiber of the projection is a single homothety equivalence class.
\end{exer}

\subsubsection{Fibred products}

\begin{fact}
(i) The category $k_H$-$An$ possesses a fibred product $Y\times_XZ$ which agrees with the fibred product in any category $k_{H'}$-$An$ for $H\subseteq H'$ and in the category of $k$-affinoid spaces. In particular, if $X=\calM(A)$, $Y=\calM(\calB)$ and $Z=\calM(\calC)$ then $\calM(\calB\wtimes_\calA\calC)\toisom Y\times_XZ$

(ii) Let $f{\colon}Y\to X$ and $g{\colon}Z\to X$ be morphisms of $k$-analytic spaces, and assume that $X=\cup_i (X_i)$, $f^{-1}(X_i)=\cup_j Y_{ij}$ and $f^{-1}(X_i)=\cup_k Z_{ik}$ are admissible coverings by affinoid domains. Then $Y\times_XZ$ admits an admissible covering by affinoid domains $Y_{ij}\times_{X_i}Z_{ik}$.
\end{fact}

Actually, the second part of this result indicates how the fibred product is constructed.

\subsubsection{The category $An$-$k$}
Often one also needs to consider morphisms between analytic spaces defined over different fields. For example, to define fibers of morphisms or ground field extensions. For this Berkovich introduces the category of analytic $k$-spaces.

\begin{defin}\label{an-k-def}
(i) An {\em analytic $k$-space} is a pair $(X,K)$ where $K$ is a non-archimedean $k$-field and $X$ is a $K$-analytic space.

(ii) A morphism $(Y,L)\to (X,K)$ consists of a bounded homomorphism $\phi:K\into L$, a continuous and $G$-continuous map $f:Y\to X$, affinoid nets (or atlases) $Y=\cup_{i\in I}Y_i$ and $X=\cup_{i\in I}X_i$ such that $f(Y_i)\subseteq X_i$, and bounded homomorphisms $\phi_i:\calO_{X_G}(X_i)\to\calO_{Y_G}(Y_i)$ that extend $\phi$, agree with $f:Y_i\to X_i$ after applying $\calM$, and pairwise agree on intersections (i.e. $\phi_i$ and $\phi_j$ agree on any $Y_k\subseteq Y_i\cap Y_j$). One identifies morphisms via refinement of atlases pretty similar to the definition of morphisms of $k$-analytic spaces.

(iii) The category of {\em analytic $k$-spaces} is denoted $An$-$k$.
\end{defin}

\subsubsection{Fibers of morphisms and base change}
\begin{defex}
(i) It follows from fact \ref{fibfact} that for any point $x\in X$ the field $\calH(x)$ is well defined: one takes any affinoid domain $x\in V$ and defines $\calH(x)$ using that domain. Obviously, $\calH(x)$ is preserved when we replace $X$ with any analytic domain containing $x$.

(ii) Assume that $f{\colon}Y\to X$ is a morphism of $k$-analytic spaces and $x\in X$ is a point. Define the fibred product $Y_x=Y\times_X\calM(\calH(x))$ as an $\calH(x)$-analytic space. Show that $Y_x=f^{-1}(x)$ set-theoretically. The space $Y_x$ is called the {\em fiber} of $f$ over $x$.

(iii) Assume that $X$ is a $k$-analytic space and $K$ is a non-archimedean $k$-field. Then one can define a universal $K$-analytic space $X_K=X\wtimes_kK$ such that any morphism $Y\to X_k$ factors uniquely through $X_K$. One simply takes an atlas $\{X_i=\calM(\calA_i)\}$ of $X$, shows that $\{X_{K,i}=\calM(\calA_i\wtimes_kK)\}$ is an atlas of a $K$-analytic space, and calls that space $X\wtimes_kK$. The construction does not depend on the atlas and provides a functor $k$-$An\to K$-$An$ called the {\em ground field extension} functor.
\end{defex}

Actually, one can unify these two constructions by saying that a fibred product $Y\times_XZ$ in $An$-$k$ exists whenever $X$ and $Y$ are $k$-analytic and $Z=\calM(K)$ with $K$ a non-archimedean field.

\begin{rem}
The definition of analytic $k$-space $X=(X,K)$ fixes a field of definition $K$. Sometimes this may be too restrictive because the "same" $X$ may admit many different fields of definition. (This is analogous to the fact that complete local rings or even non-reduced affine varieties usually admit different fields of definitions.) For example, $K_r\{T\}$ admits many bounded $k$-automorphisms $\phi$ that do not preserve $K_r$, and one may wish to view $\calM(\phi)$ as an automorphism of the analytic $k$-space $\calM(K_r\{T\})$. There is a simple way to extend the category $An$-$k$ accordingly: in Definition \ref{an-k-def}, one omits a morphism between fields of definitions and allows all morphisms $Y_i\to X_i$ corresponding to a $k$-bounded homomorphism $\calO_{X_G}(X_i)\to\calO_{Y_G}(Y_i)$.
\end{rem}

\subsection{Basic classes of analytic spaces and morphisms}

\subsubsection{Good spaces}

\begin{defin}
(i) An analytic space $X$ is called {\em good at} a point $x$ if $x$ has an affinoid neighborhood. The space $X$ is {\em good} if it is good at all its points.

(ii) The sheaf $\calO_X$ is defined as the restriction of $\calO_{X_H}$ to the usual topology of $X$.

(iii) For any point $x\in X$ we define $\kappa(x)$ as the residue field of the local ring $\calO_{X,x}$. Note that $\kappa(x)$ may change when we replace $X$ with an analytic domain containing $x$ (and this happens already in the affinoid case).
\end{defin}

Good spaces are often more convenient to work with. In particular, the definition of $\calO_X$ applies to any $X$, but it can play the role of a structure sheaf only for the class of good spaces (see Example \ref{badex}). That is why we will only consider $\calO_X$ and $\kappa(x)$ on good spaces.

\begin{fact}
For a good $k_H$-analytic space $X$ the category of coherent $\calO_{X_H}$-modules is equivalent to the category of {\em coherent} $\calO_X$-modules (i.e. $\calO_X$-modules locally isomorphic to a quotient of $\calO_X^n$).
\end{fact}

\begin{exex}\label{badex}
(i) If an analytic space $X$ is good at a point $x\in X$ then $\calH(x)$ is the completion of $\kappa(x)$.

(ii) Let $X$ be a closed unit disc with doubled open unit disc as in Exercise \ref{glueex}, and let $x$ be its maximal point. Show that in some but not all cases $k\toisom\calO_{X,x}$. Thus, the usual topology is too crude to allow non-constant functions in a neighborhood of $x$. Clearly, $X$ is not good and $\wh{\kappa(x)}\subsetneq\calH(x)$ in this case.
\end{exex}

Here are simplest separated examples of non-good spaces. They will be studied further in \S\ref{RZsec} (and, at least, one has to use Raynaud's theory from \S\ref{raysec} to show that these examples are non-good spaces.)

\begin{exam}\label{nogoodexam}
(i) Assume that $\ur$ is a tuple of $n>1$ positive numbers. A closed polydisc $E(0,\ur)$ of polyradius $\ur$ with removed open polydisc of polyradius $\ur$ is a compact not good analytic domain in $X\subset E(0,\ur)$.

(ii) In the affine plane $X=\MSpec(k[S,T])$ consider the affinoid domains $V_1=X\{r\le|S|\le 1,|ST|\le1\}$ and $V_2=X\{1\le |S|\le r^{-1},|S^{-1}T|\le 1\}$ with $0<r<1$. One can show that the compact analytic domain $V=V_1\cup V_2$ is not good at the maximal point of the unit polydisc $X\{S,T\}$. Moreover, if we consider the natural projection $f:V\to\bfA^1_k=\MSpec(k[S])$ then $V_1$ and $V_2$ are the preimages in $V$ of the affinoid domains $A(0;r,1)$ and $A(0;1,r^{-1})$. Since the union of these two annuli is an affinoid space (it is $A(0;r,r^{-1}))$), we see that no reasonable notion of good or affinoid morphisms exists.
\end{exam}

\subsubsection{Finite morphisms, closed immersions and Zariski topology}

\begin{defin}
A morphism $f{\colon}Y\to X$ is called a {\em closed immersion} (resp. {\em finite}) if there exists an admissible covering by affinoid domains $X=\cup_i X_i$ such that $Y_i=X_i\times_XY$ are $k$-affinoid and the homomorphisms of Banach algebras $\calO_{X_G}(X_i)\to\calO_{Y_G}(Y_i)$ are surjective (resp. finite) and admissible.
\end{defin}

We have already discussed in Exercise \ref{admex} why admissibility is essential when the valuation on $k$ is trivial.

\begin{exer}\label{finex}
(i) Assume that $f\colon Y\to X$ is finite. Then for any affinoid domain $\calM(\calA)=X'\subseteq X$ its preimage $Y'=Y\times_XX'$ is affinoid, say $Y'=\calM(\calB)$, and the homomorphism $\calA\to\calB$ is surjective (resp. finite) and admissible. (Hint: $\calO_{Y'_G}$ is a coherent $\calO_{X'_G}$-algebra.)

(ii) The class of closed immersions (resp. finite morphisms) is closed under compositions, base changes and ground field extensions.
\end{exer}

\begin{defin}
Any subset $Z\subseteq X$ that is the image of a closed immersion is called {\em Zariski closed}. The complement of such set is called {\em Zariski open}.
\end{defin}

\begin{exer}
A point $x\in X$ is Zariski closed if and only if $[\calH(x):k]<\infty$. (Zariski closed points of $X$ are precisely its classical rigid points. The set of all such points is denoted $X_0$.)
\end{exer}

When working with Zariski topology one must be very careful because it becomes stronger when passing to analytic domains (even open ones). In other words, coherent ideals on an open subspace do not have to extend to the whole space. Such phenomenon does not occur for algebraic varieties (for obvious reasons) but does occur for formal varieties.

\begin{exam}\label{badbound}
(i) Give an example of a $k$-analytic space $X$ with an open subspace $U$ and a closed subspace $Z\subset U$ which does not extend to the whole $X$. (Hint: take $X=\calM(k\{T,S\})$ the unit polydisc, $U$ an open polydisc of polyradius $(1,r)$ with $r<1$ and $Z$ given by $T-f(S)$, where $f(S)$ has radius of convergence between $r$ and $1$.)

(ii) An example of Ducros. Fix $r\notin\sqrt{|k^\times|}$ with $0<r<1$. Consider a polydisc $X=\calM(k\{T,S\})$ with an affinoid domain $V=\calM(k\{r^{-1}T,rT^{-1},S\})$, which is a unit $K_r$-disc and a non-strict $k$-surface (the product of the unit $k$-disc with the irrational $k$-annulus $\calM(K_r)$). Using the hint from (i) find a Zariski closed point $x\in V$ with $\calH(x)\toisom K_r$ such that the ideal of $x$ in $V$ does not extend to any neighborhood of $x$ in $X$. (In a sense, $x$ is a $k$-curve in the $k$-surface $X$ which cannot be extended, so $x$ is Zariski closed only in a sufficiently small domain in $X$.) Show that in this case
$\calO_{X,x}$ is a dense subfield of $K_r$ and hence the character $\chi_{X,x}{\colon}k\{T,S\}\to\calH(x)$ is injective and hence flat. On the other hand, its base change with respect to the homomorphism $k\{T,S\}\to k\{r^{-1}T,rT^{-1},S\}$ is not flat because the character $\chi_{V,x}$ has a non-trivial kernel. In particular, one cannot define a reasonable class of flat morphisms between $k$-affinoid spaces just by saying that $\calM(f)$ is flat whenever $f$ is flat. (One can show that this approach works well for strictly $k$-affinoid spaces; see also Remark \ref{strrem}(ii).)
\end{exam}

\subsubsection{Separated morphisms}
\begin{defin}
(i) A morphism $f{\colon}Y\to X$ is {\em separated} if the diagonal morphism $\Delta_{Y/X}{\colon}Y\to Y\times_XY$ is a closed immersion.

(ii) A $k$-analytic space is {\em separated} if so is its morphism to $\calM(k)$.
\end{defin}

\begin{exer}
(i) Formulate and prove the basic properties of separated morphisms analogous to the properties of separated morphisms of schemes. In particular, show that in a separated $k$-analytic space $X$ the intersection of two $k$-affinoid domains is $k$-affinoid.

(ii) Prove Fact \ref{faithfact} for the subcategories of separated objects in $k_H$-$An$ and $k$-$An$. (Hint: if $X$ and $Y$ are $H$-strict then any closed subspace in $X\times Y$ is affinoid by Exercise \ref{finex} and hence $H$-strict.)
\end{exer}

Non-separatedness of a space can be of two sorts, as is illustrated by the following example.

\begin{exex}
(i) Let $X$ be the closed unit disc with doubled open unit disc. Show that $X$ is not locally separated at its maximal point $x$ (i.e. any $k$-analytic domain which is a neighborhood of $x$ is not separated). In particular, $X$ is a non-good $k$-analytic Hausdorff space.

(ii) Show that the closed unit disc $Y$ with doubled origin is not separated but is locally separated at all its points. Moreover, $Y$ is a good non-Hausdorff $k$-analytic space.
\end{exex}

\subsubsection{Boundary and proper morphisms}\label{propsec}
\begin{defin}
(i) The {\em relative interior} $\Int(Y/X)$ of a morphism $f{\colon}Y\to X$ is the set of all points $y\in Y$ such that for any affinoid domain $U\subseteq X$ containing $x=f(y)$ there exists an affinoid domain $V\subseteq f^{-1}(U)$ such that $V$ is a neighborhood of $x$ in $f^{-1}(U)$ and $x\in\Int(V/U)$. The complement $\partial(Y/X)=Y\setminus\Int(Y/X)$ is called the {\em relative boundary} and we say that $f$ is {\em boundaryless} or without boundary if $\partial(Y/X)$ is empty (in \cite{berihes} such morphisms are called
"closed").

(ii) A morphism $f{\colon}Y\to X$ is {\em proper} if it is boundaryless and compact (i.e. the preimage of a compact domain is compact).
\end{defin}

As usual, the absolute analogs of these notions are defined relatively to $\calM(k)$. Note that (part (i) of) this definition agrees with our earlier definitions from \S\ref{boundsec}.

\begin{exex}
(i) A $k$-analytic space has no boundary if and only if any its point $x$ possesses an affinoid neighborhood $U$ such that $x\in\Int(U)$. For example, an open polydisc, $\bfP^n_k$ and $\bfA^n_k$ have no boundary, and a closed polydisc has a boundary. So far, $\bfP^n_k$ is the only example of a proper non-discrete $k$-analytic space we have considered. It follows that any {\em projective $k$-analytic space}, i.e. a closed subspace of $\bfP^n_k$ is also proper.

(ii) Any boundaryless $k$-analytic space $X$ is good. If the valuation on $k$ is non-trivial then $X$ is also strict.

(iii) A morphism between affinoid spaces is proper if and only if it is finite. Any finite morphism is proper.

(iv) A boundaryless morphism is separated if and only if the preimage of any Hausdorff domain is Hausdorff. In particular, proper morphisms are separated.

(v) There is a theory of analytic tori which is parallel in some part to the classical theory of complex tori. An {\em analytic torus} is defined as the quotient $T_\Lam=\bfG_m^d/\Lam$ where $\Lam=\oplus_{i=1}^d\lam_i^\bfZ$ is a multiplicative lattice such that $|\Lam|$ is a lattice in $(\bfR_+^\times)^d$. It is easy to see that $T_\Lam$ is a proper analytic space. There also is a an analog of Riemann's positivity conditions on $\Lam$ that are necessary and sufficient for $T_\Lam$ to be algebraic (and even projective). Similarly to the complex case, if $d=1$ then $T_\Lam$ is always algebraic, but a generic two dimensional torus is not algebraic.
\end{exex}

\begin{fact}\label{boundfact}
(i) If $Y$ is an analytic domain in $X$ then $\Int(Y/X)$ is the topological interior of $Y$ in $X$.

(ii) Boundaries are $G$-local on the base, i.e. given an admissible covering of $X$ by affinoid domains $X_i$ one has that $\partial(Y/X)=\cup_i\partial(X_i\times_XY/X_i)$.

(iii) The classes of proper morphisms and morphisms without boundary are $G$-local on the base and are preserved by compositions, base changes and ground field extensions.

(iv) If $f{\colon}Y\to X$ is a separated boundaryless morphism and $X$ is $k$-affinoid then for any affinoid domain $U\subseteq Y$ there exists a larger affinoid domain $V\subseteq Y$ such that $U\subseteq\Int(V/X)$ and $U$ is a Weierstrass domain in $V$.
\end{fact}

\begin{rem}\label{proprem}
Surprisingly enough, already (ii) is really difficult. It turns out that when one wants to show that various morphisms have no boundary, the difficult part of the proof is to show that the preimage of an affinoid domain under these morphisms is a good domain. Once this is established, one can use the theory of boundaries for affinoid spaces as outlined in \S\ref{boundsec}. See also Remark \ref{kirem} below.
\end{rem}

Similarly to algebraic and complex analytic geometries, coherence is preserved by higher direct images with respect to proper morphisms.

\begin{fact}[Kiehl's theorem on direct images]\label{kith}
If $f{\colon}Y\to X$ is a proper morphism between $k$-analytic spaces and $\calF$ is a coherent
$\calO_{Y_G}$-modules then the $\calO_{X_G}$-modules $R^if_*(\calF)$ are coherent.
\end{fact}

Note that we use here that $f$ is a compact map because otherwise $f_*(\calF)$ is not a sheaf of Banach $\calO_{X_G}$-modules.

\begin{rem}\label{kirem}
Kiehl introduced the notion of proper morphisms and proved the above result (for rigid spaces) in \cite{kiehl}. One can easily show that our definition of proper morphisms (in the strict case) is equivalent to the original Kiehl's definition. The definition of proper morphisms is designed so that the theorem on direct images can be proved rather easily and naturally (one computes \v{C}ech complexes and shows that certain differentials are compact operators). As was already remarked, it is very difficult to establish some other properties, that one might expect to be more foundational. For example, the fact that proper morphisms are preserved by compositions was open for more than twenty years (for a discretely valued $k$ this was proved in \cite{lut} and the general case was established in \cite{temred1} and \cite{temred2}).
\end{rem}

\subsubsection{Smooth and \'etale morphisms}
\begin{defin}
(i) A finite morphism $f{\colon}Y\to X$ is {\em \'etale} if for any affinoid domain $U\subseteq X$ and its preimage $V=f^{-1}(U)$ the finite homomorphism of $k$-affinoid algebras $\calO_{X_G}(U)\to\calO_{Y_G}(V)$ is \'etale. (Recall that $V$ is affinoid and finite over $U$ by Example \ref{finex}(i).)

(ii) In general, a morphism $f{\colon}Y\to X$ is {\em \'etale} if locally (on $Y$) it is finite \'etale. Namely, for any point $y\in Y$ there exist neighborhoods $V$ of $y$ and $U$ of $f(y)$ such that $f$ restricts to a finite \'etale morphism $V\to U$.

(iii) A morphism $f{\colon}Y\to X$ is {\em smooth} if it can be represented as an \'etale morphism $Y\to\bfA^n_X$ followed by the projection.
\end{defin}

\begin{exer}
(i) Any smooth morphism (e.g. an \'etale morphism) is boundaryless.

(ii) If $V$ is an analytic domain in $X$ then the embedding $V\into X$ is \'etale if and only if it is an open immersion.

(iii)* Any smooth morphism is an open map.
\end{exer}

\begin{fact}
The classes of \'etale and smooth morphisms are closed under compositions, base changes and ground field extensions. Also, it follows from Kiehl's theorem that \'etaleness is $G$-local on the base. Probably, smoothness is not $G$-local on the base.
\end{fact}

This definition of \'etale and smooth morphisms is analogous to a complex analytic definition but it does not apply to nice morphisms with boundaries. For example, the closed unit disc is not smooth at its maximal point. The following definition is a natural generalization to the case when there are boundaries. We give it for the sake of completeness, but do not discuss all results one should prove to show that it really makes sense.

\begin{defin}
(i) A morphism $f{\colon}Y\to X$ between strictly $k$-analytic spaces is {\em rig-smooth} if the restriction of $f$ on $\Int(Y/X)$ is smooth.

(ii) In general, a morphism $f{\colon}Y\to X$ is {\em rig-smooth} if so is some (and then any) ground field extension $f_\ur:=f\wtimes_kK_\ur$ such that $Y_\ur$ and $X_\ur$ are strictly $K_r$-analytic.

(iii) A rig-smooth morphism with discrete fibers is called {\em quasi-\'etale}.
\end{defin}

\begin{rem}\label{strrem}
(i) Alternatively, one can define quasi-\'etale morphisms directly and then rig-smooth morphisms are the morphisms that locally split into the composition of a quasi-\'etale morphism with the projection $\bfA^n_X\to X$.

(ii) We have to extend the ground field in the general case because $\Int(Y/X)$ can be too small to test (any sort of) smoothness. A good example of such situation was studied in Exercise \ref{badbound}(ii).

(iii) The same problem as in (ii) happens when one wants to introduce flatness. A reasonable theory of flatness was developed very recently by Ducros. In the strict case one gives a naive definition, and in general $f$ is called {\em flat} if so is its strictly analytic ground field extension.

(iv) Ducros also proves that $f$ is rig-smooth if and only if it is flat, the coherent $\calO_{Y_G}$-module of continuous differentials $\Omega^1_{Y_G/X_G}$ is locally free, and locally the rank of $\Omega^1_{Y_G/X_G}$ equals to the relative dimension. (The sheaf $\Omega^1_{Y_G/X_G}$ admits the following local description: if $X=\calM(\calA)$ and $Y=\calM(\calB)$ then $\Omega^1_{Y_G/X_G}$ corresponds to the module $I/I^2$ where $I=\Ker(\calB\wtimes_\calA\calB\to\calB)$).
\end{rem}

\subsection{Topological properties}

\subsubsection{Basic properties}

\begin{fact}
(i) Any connected $k$-analytic space is pathwise connected.

(ii) Any point has a fundamental family of neighborhoods which are compact and pathwise connected analytic domains.

(iii) The topological dimension of a paracompact $X$ is at most $\dim(X)$ and both are equal in the strict case.
\end{fact}

This fact was proved in \cite{berbook} (and the argument is correct, although by a misunderstanding some mathematicians thought that part (i) was not proved). Let us say few words about the proof of (i). We have already checked this fact for an affine line and one easily deduces the case of a polydisc. Studying finite covers of polydiscs one obtains the case of strictly $k$-analytic spaces. The general case is deduced by descent from an appropriate $X\wtimes_k K_\ur$. Next let us consider examples of analytic spaces with bad topological properties.

\begin{exam}\label{badexam}
(i) Assume that $\tilk$ is uncountable (e.g. $k=\bfC((T))$). Let $U$ be the open subset of a closed two-dimensional polydisc $E$ obtained by removing the maximal point of $E$. Then $U$ is not paracompact, i.e. it possesses open covers that do not admit locally finite refinements.

(ii) There exist (see Exercise \ref{infex}) examples of $k$-analytic curves $C$ such that $C$ is a double covering of an open unit disc and $C$ is a closed subspace in a two-dimensional open unit polydisc, but the first Betti number of $C$ is infinite. In particular, one can construct $C$ so that it can be retracted onto one of its subset $\Delta$, which is a graph with infinitely many loops.
\end{exam}

\subsubsection{Contractions}
In some cases one can construct by hands a retraction of a $k$-analytic space $X$ onto a subset $S$ which is of topologically finite type. One such method is to find an action of a $k$-affinoid group $G$ on $X$ with a continuous family of affinoid subgroups $\{G_t\}_{t\in[0,1]}$ such that $G_0=\{e\}$, $G_1=G$ and for each point $x\in X$ each orbit $G_tx$ is affinoid and possesses exactly one maximal point $x_t$. Then $(x,t)\mapsto x_t$
defines a deformational retraction of $X$ onto some its subset, which is very small in some examples. Two good examples of such $G$ are as follows: a closed unit polydisc $\calM(k\{\uT\})$ with an additive group structure, and a product of unit annuli $G_{m,1}^n=\calM(k\{T_1,T_1^{-1}\.T_n,T_n^{-1}\})$ with the multiplicative group structure. The
groups $G_t$ with $t<1$ are the polydiscs of polyradius $(t\. t)$ with center at $0$ or $1$, respectively. The action of the torus is much more important because tori play important role in the theory of reductive groups (see \cite[\S5]{berbook} for the connection to Bruhat-Tits buildings). So, we give the most fundamental example of a contraction by a torus action.

\begin{exex}
(i) Show that $\bfR_+^n$ embeds into $\bfA^n_k$ so that $\ur=(r_1\. r_n)$ goes to the semivaluation $\|\ \|_\ur$ (the maximal point of $E(0,\ur)$).

(ii)* Show that the action of $G_{m,1}^n$ on $\bfA^n_k$ contracts it onto $\bfR_+^n$. Moreover, the retraction can be explicitly described by the formula $|f(x_t)|=\max_{i\in\bfN^n}|\partial_if(x)|t^i$ where $f(\uT)\in k[T_1\. T_n]$ and $\partial_i{\colon}k[\uT]\to k[\uT]$ for $i\in\bfN^n$ is the logarithmic differential operator
$\frac{\uT^i}{i!}\frac{\partial^i}{d\uT^i}$.
\end{exex}

Note also that Berkovich devoted a separate paper \cite{bercontr} to proving the following very difficult result.

\begin{fact}\label{contrfact}
Any analytic domain in a smooth $k$-analytic space is locally contractible.
\end{fact}

\subsubsection{Topological type of analytic spaces}\label{topsec}
We saw in \S\ref{aflinesec} that the affine line is a sort of an infinite tree. The topological structure of general analytic spaces is much more complicated and somewhat mysterious. A major progress in its understanding was achieved very recently by Hrushovski-Loeser in \cite{HL} via model-theoretic methods, and here is one of the main applications of their theory to analytic spaces.

\begin{fact}\label{topconj}
Let $\oX$ be a projective $n$-dimensional $k$-analytic space, and assume that $X=V\cap Y$ is the intersection of a Zariski open subspace $V\subset\oX$ with a compact analytic domain $Y\subset\oX$. Then $X$ contains a family of topological spaces $\{S_i\}_{i\in I}$ filtered by inclusion such that each $S_i$ is homeomorphic to a finite simplicial complex of dimension at most $n$ and there exists a projective family of maps $f_{ij}{\colon}S_i\to S_j$ (for each pair $S_j\subseteq S_i$) such that $X\toisom\projlim_{i\in I} S_i$. Moreover, this family extends
to a compatible family of deformational retractions $\Phi_{ij}{\colon}S_i\times[0,1]\to S_j$ with $\Phi_{ij}(x,1)=f_{ij}(x)$ which induce deformational retractions of $X$ on each of $S_i$'s.
\end{fact}

Since any rig-smooth analytic space can be locally embedded as a subdomain in a projective variety, this result immediately implies Fact \ref{contrfact}. Moreover, this result is of global nature and it treats most types of singularities as well (although, there exist non-algebraizable singularities that cannot be locally embedded into varieties). It seems natural to expect that the projective variety $\oX$ in Fact \ref{topconj} can be replaced with an arbitrary compact $k$-analytic space, but in such generality this is a widely open conjecture. (Note that in order to exclude bad spaces discussed in Example \ref{badexam}, some compactness assumption should be present in the formulation.)

\section{Relation to other categories}\label{5}
This section contains various material, and some of its subsections are rather advanced. We place it before section \S\ref{cursec} on analytic curves because some results of sections \S\S\ref{ansec}--\ref{raysec} will be used to study curves. So, the reader can look through the first three subsections and go directly to \S\ref{cursec}.

\subsection{Analytification of algebraic $k$-varieties}\label{ansec}

\subsubsection{The analytification functor}

Let $k$-$Var$ be the category of algebraic $k$-varieties (i.e. schemes of finite type over $k$). We are going to describe a construction of an analytification functor $k$-$Var\to k$-$An$. The analytification of a morphism $f{\colon}\calY\to\calX$ will be denoted $f^\an{\colon}\calY^\an\to\calX^\an$. For $\calX=\Spec(k[T_1\. T_n])$ we set $\calX^\an=\bfA^n_k=\MSpec(k[\uT])$. For any quotient $A=k[\uT]/I$ the analytification of $\calY=\Spec(A)$ is the closed subspace of $\calX^\an$ defined by vanishing of $I\calO_{\calX^\an}$.

\begin{exer}
(i) Prove that this definition is independent of choices. Also, show that $\calY^\an=\MSpec(A)$ is the set of all real semivaluations on $A$ bounded on $k$. (Hint: two embeddings of $\calY$ into affine spaces are dominated by a third such embedding.)

(ii) Extend this to a functor from the category of affine $k$-varieties to the category of boundaryless $k$-analytic spaces.

(iii) Show that the latter functor takes open immersions to open immersions and hence extends (by gluing) to an analytification functor $k$-$Var\to k$-$An$. Show that any analytification is a boundaryless space.

(iv) Show that $(\Proj(A))^\an\toisom\MProj(A)$ for any graded finitely generated $k$-algebra $A$. In particular, $\MProj(A)$ is a projective analytic variety.
\end{exer}

\begin{fact} The analytification functor can be described via the following universal property. For any good $k$-analytic space $Y$ let $F_\calX(Y)$ be the set of morphisms of locally ringed spaces $(Y,\calO_Y)\to(\calX,\calO_\calX)$. Then $X=\calX^\an$ is the $k$-analytic space that represents $F_\calX$.
\end{fact}

In particular, a morphism $\pi_\calX{\colon}(\calX^\an,\calO_{\calX^\an})\to(\calX,\calO_\calX)$ arises.

\begin{fact}
$\calX^\an(K)\toisom\calX(K)$ for any non-archimedean $k$-field $K$, in particular, $\pi_\calX$ is surjective.
\end{fact}

\begin{defin}
For a coherent $\calO_\calX$-module $\calF$ the module $\calF^\an=\pi_\calX^*(\calF)$ is a coherent $\calO_{\calX^\an}$-module called the {\em analytification} of $\calF$.
\end{defin}

\subsubsection{GAGA}
The analytification functor preserves almost all properties of varieties and their morphisms, and here is a
(partial) list.

\begin{fact}
Let $f{\colon}\calY\to\calX$ be a morphism between algebraic $k$-varieties. Then $f$ satisfies one of the following properties if and only if so does $f^\an$: smooth, \'etale, finite, closed immersion, open immersion, isomorphism, proper, separated.
\end{fact}

\begin{fact}\label{gaga}
(i) For a proper variety $\calX$ the analytification functor induces an equivalence
$\Coh(\calO_\calX)\toisom\Coh(\calO_{\calX^\an})$.

(ii) The functor $\calX\mapsto\calX^\an$ is fully faithful on the category of proper varieties
\end{fact}

\begin{exer}
(i) The assertions of Fact \ref{gaga} do not hold for general algebraic varieties.

(ii) For a proper variety $\calX$, the analytification functor induces an equivalence between the categories of finite (resp. finite \'etale) $\calX$-schemes and $\calX^\an$-spaces.

(iii) Any projective $k$-analytic space $X$ is algebraizable by a projective $k$-variety $\calX$ (i.e. $X\toisom\calX^\an$).
\end{exer}

When the valuation on $k$ is trivial, the properness assumption can be eliminated. Let $X_t\subset X$ be the set of points $x\in X$ with trivially valued completed residue field $\calH(x)$.

\begin{fact}
Assume that the valuation on $k$ is trivial.

(i) $\calX_t^\an\toisom\calX$.

(ii) The analytification functor is fully faithful.

(iii) For a variety $\calX$ the analytification functor induces an equivalence of categories
$\Coh(\calO_\calX)\toisom\Coh(\calO_{\calX^\an})$.
\end{fact}

\subsection{Generic fibers of formal $\kcirc$-schemes}

\subsubsection{Reminds on formal schemes}
\begin{defin}
Let $A$ be a ring with an ideal $I$.

(i) The {\em $I$-adic topology} on $A$ is generated by the cosets $a+I^n$.

(ii) The {\em separated $I$-adic completion} is defined as $\hatA=\projlim_n A/I^n$.

(iii) $A$ is {\em $I$-adic} if $A\toisom\hatA$.

(iv) Any ideal $J$ with $I^n\subseteq J$ and $J^n\subseteq I$ for large enough $n$ is called {\em ideal of definition} of $A$. It generates the same topology and can be used instead of $I$ in all definitions.
\end{defin}

\begin{exex}
(i) If $\kcirc$ is a real valuation ring with fraction field $k$ and $\pi\in\kcirccirc$ is any non-zero element then the $(\pi)$-adic completion of $\kcirc$ is the ring of integers $\hatkcirc$ of the completion of $k$.

(ii) The separated $(\kcirccirc)$-completion of $k$ is either $\hatkcirc$ or $\tilk$. Moreover, the first possibility occurs only when $k$ is discrete or trivially valued.
\end{exex}

\begin{defin}
(i) The {\em formal spectrum} $\gtX=\Spf(A)$ of an $I$-adic ring $A$ is the set of open prime ideals of $A$ with the topology generated by the sets $D(f)$, where $D(f)$ is the non-vanishing locus of an element $f\in A$.

(ii) Each $D(f)$ is homeomorphic to $\Spf(A_{\{f\}})$, where the {\em formal localization} $A_{\{f\}}$ is the universal $I$-adic $A$-algebra with inverted $f$.

(iii) The structure sheaf $\calO_\gtX$ is the sheaf of topological rings determined by the condition $\calO_\gtX(D(f))=A_{\{f\}}$. The topologically ringed space $(\gtX,\calO_\gtX)$ is called the {\em affine formal scheme} associated with $A$.

(iv) The {\em closed fiber} (or special fiber) $\gtX_s$ of $\gtX$ is the reduction of $\Spec(A/J)$ for any ideal of definition $J$.
\end{defin}

\begin{exer}\label{formex}
(i) Show that $A_{\{f\}}\toisom\projlim_n(A/I^n)_f$ and $A\{T\}/(Tf-1)\toisom A_{\{f\}}$, where
$A\{T\}=\projlim_n(A/I^n)[T]$ is the ring of convergent power series over $A$.

(ii) If $\pi^n\in I$ for some $n$ then $A_{\{\pi\}}=0$. Thus, formal localization at topologically nilpotent element has the same effect as inverting a nilpotent element in a ring.

(iii) Show that $\gtX_s$ does not depend on the ideal of definition and $|\gtX_s|\toisom|\gtX|$. Actually, $\gtX$ can be viewed as the inductive limit of schemes $(\gtX_s,\calO_\gtX/J)$ where $J$ runs through the ideals of definition.
\end{exer}

\begin{defin}
A general {\em formal scheme} is a topologically ringed space $(\gtX,\calO_\gtX)$ which is locally isomorphic to affine formal schemes. Morphisms of such creatures are morphisms of topologically ringed spaces that induce local homomorphisms on the ring-theoretical stalks. The notions of ideals of definitions and of the closed fiber are extended to the general formal schemes in the obvious way.
\end{defin}

\begin{defex}
(i) The $n$-dimensional affine space over an adic ring $A$ is defined as $\bfA^n_A=\Spf(A\{T_1\. T_n\})$.

(ii) A formal scheme over an $I$-adic ring $A$ is of {\em (topologically) finite presentation} (resp. {\em special}) if it is locally of the form $\Spf(A\{T_1\. T_n\}/(f_1\. f_m))$ (resp. $\Spf(A\{T_1\. T_n\}[[S_1\. S_l]]/(f_1\. f_m))$).
\end{defex}

\subsubsection{Generic fibers of formal $\kcirc$-schemes of finite type}
In this section we are going to define a {\em generic fiber functor} $\eta$ which assigns to a formal $\kcirc$-scheme $\gtX$ of locally finite type a Hausdorff strictly $k$-analytic space $\gtX_\eta$ (even when the valuation is trivial). Intuitively, $\gtX_\eta$ is the "missing generic fiber of $\gtX$" and when $k$ is non-trivially valued it is defined by inverting a non-zero element $\pi\in\kcirccirc$. (By Exercise \ref{formex}(ii) we kill any formal $\kcirc$-scheme by such an operation, so it is not surprising that $\gtX_\eta$ is not a formal scheme but leaves in another category.) If $k$ is trivially valued then we set $\pi=0$ to uniformize the exposition.

In general, the definition of $\eta$ is very similar to the definition of the analytification. One defines $(\bfA^n_{\kcirc})_\eta$ to be the closed unit polydisc $E^n(0,1)$. An affine scheme given by vanishing of $f_1\. f_n$ in $\bfA^n_{\kcirc}$ is defined as the closed subspace in $E^n(1,0)$ given by vanishing of $f_i$'s. In general, $\eta$ is defined via gluing. Let us realize this program with some details.

\begin{defin}
If $A=\kcirc\{T_1\. T_n\}/I$ is a $\pi$-adic ring with a finitely generated $I$ then
$\calA=A_\pi=A\otimes_{\kcirc} k$ is a $k$-affinoid algebra isomorphic to $k\{\uT\}/Ik\{\uT\}$. For the affine formal scheme $\gtX=\Spf(A)$ we set $\gtX_\eta=\calM(\calA)$.
\end{defin}

\begin{exer}\label{torex}
Assume that $A$ has no $\pi$-torsion, and so $A$ embeds into $\calA$.

(i) The integral closure of $A$ in $\calA$ is $\calAcirc$.

(ii) If the valuation on $k$ is non-trivial and $A$ is reduced then $A$ is the unit ball for a Banach norm $|\ |$ on $\calA$. This means that $|\ |$ is equivalent to $\rho_\calA$ or, equivalently, $\pi^n\calAcirc\subseteq A\subseteq\calAcirc$ for large enough $n$. (Hint: use Fact \ref{basicaffact}(iii).)

(iii) Formal localization is compatible with inverting $\pi$. Namely, $(A_{\{f\}})_\pi\toisom A_\pi\{f^{-1}\}$. (Hint: use Fact \ref{basicaffact}.)
\end{exer}

Part (iii) of the above exercise implies that $\eta$ (defined for affine formal schemes) takes open immersions to embeddings of affinoid domains. Now we can define the functor $\eta$ in general.

\begin{defex}
(i) If a separated formal scheme $\gtX$ of finite type over $\kcirc$ is glued from open subschemes $\gtX_i$ along the intersections $\gtX_{ij}$ then the gluing of $(\gtX_i)_\eta$ along $(\gtX_{ij})_\eta$ is possible by Exercise \ref{glueex}(ii) and the obtained $k$-analytic space is set to be $\gtX_\eta$.

(ii) For a general formal scheme $\gtX$ we repeat the same construction but with $\gtX_{ij}$ being separated formal schemes now.

(iii) Check that this construction defines the promised generic fiber functor (in particular, it extends to morphisms).
\end{defex}

As one might expect, $\eta$ preserves (or naturally modifies) almost all properties of morphisms.

\begin{fact}
Let $f{\colon}\gtY\to\gtX$ be a morphism between formal $\kcirc$-schemes without $\pi$-torsion. If $f$ is an isomorphism, separated, proper, a closed immersion, finite \'etale, then $f_\eta$ is so. If $f$ is an open immersion, \'etale, or smooth, then $f_\eta$ is a compact analytic domain embedding, quasi-\'etale, or rig-smooth, respectively.
\end{fact}

Excluding preservation of properness, all these claims are simple. The remaining claim is really difficult, though this is not so surprising in view of other problems with properness discussed in \S\ref{propsec}. Actually, to prove this claim is essentially equivalent to prove the other difficult properties of properness listed in \S\ref{propsec}.

As for the opposite implications, at first glance, one cannot expect that something can be proved in that direction. For example, a generically finite morphism does not have to be finite, etc. However, the following result holds true and its proof is relatively simple.

\begin{fact}
If $f_\eta$ is proper or separated then so is $f$.
\end{fact}

Finally, there exists an anti-continuous {\em reduction map} $\pi_\gtX{\colon}\gtX_\eta\to\gtX_s$ defined similarly to the affinoid reduction map.

\begin{defex}
Check that for any point $x\in\gtX_\eta$ its character $\chi_x{\colon}\calM(\calH(x))\to\gtX_\eta$ is induced by a morphism $\chicirc_x{\colon}\Spf(\calH(x)^\circ)\to\gtX$. This induces a point $\tilchi_x{\colon}\Spec(\wHx)\to\gtX_s$ on the closed fiber and hence gives rise to a map $\pi_\gtX$.
\end{defex}

\begin{rem}
The reduction map is an analog of the following specialization construction. If $X$ is a scheme over a henselian valuation ring (e.g. $\kcirc$), $(X_\eta)_0$ is the set of closed points of the generic fiber and $(X_s)_0$ is the set of closed points of the closed fiber then specialization induces a map $(X_\eta)_0\to(X_s)_0$.
\end{rem}

\subsubsection{Relation to the analytification}
If $k$ is trivially valued then the analytification and the generic fiber constructions provide two functors from the category of $k$-varieties to the category of $k$-analytic spaces. More generally, for any $k$ we have two functors $\calF$ and $\calG$ from the category of $\kcirc$-schemes of finite type to the category of $k$-analytic spaces: $\calG(X)=(X_\eta)^\an$ and $\calF(X)=(\hatX)_\eta$. In the first case, we first pass to the generic fiber of the morphism $X\to\Spec(\kcirc)$ and then analytify the obtained $k$-variety. In the second case, we first complete $X$ and then take the generic fiber of the obtained formal $\kcirc$-scheme of finite type.

\begin{exer}
(i) Assume that $X=\Spec(A)$ and $f_1\. f_n$ generate $A$ over $\kcirc$. Show that $(\hatX)_\eta$ can be naturally identified with the affinoid domain in $(X_\eta)^\an$ defined by the conditions $|f_i|\le 1$.

(ii) Let $\calF^{\rm sep}$ and $\calG^{\rm sep}$ be the restrictions of the functors $\calF$ and $\calG$ onto the category of separated $\kcirc$-schemes of finite type. Extend the construction of (i) to a morphism of functors $\phi{\colon}\calF^{\rm sep}\to\calG^{\rm sep}$ which is a compact embedding of a strictly analytic domain (i.e. each morphism $\phi(X){\colon}\calF(X)\to\calG(X)$ is embedding of a compact analytic domain).

(iii) Show that when restricted to proper $\kcirc$-schemes $\phi$ induces an isomorphism of functors, i.e. the embedding of the analytic domain $\phi(X){\colon}(\hatX)_\eta\into(X_\eta)^\an$ is an isomorphism for a $\kcirc$-proper $X$.

(iv) Show that $\phi$ extends to non-separated $\kcirc$ schemes but then it does not have to be embedding of an analytic domain. (Hint: if $X$ is the relative affine line over $\kcirc$ with doubled origin then $\calF(X)$ is not locally separated and hence cannot be embedded into the good (although not Hausdorff) space $\calG(X)$.)
\end{exer}

\subsubsection{Generic fibers of $\kcirc$-special formal schemes}
For completeness, we discuss briefly how the generic fiber functor extends to all $\kcirc$-special formal schemes in the case of a discretely valued (or trivially valued) ground field $k$. (The non-discretely valued case was not studied in the literature because the rings $\kcirc[[T_1\. T_n]]$ are rather pathological, e.g. they contain non-closed ideals.)

The general idea of defining $\gtX_\eta$ is actually the same: for $$\gtX=\Spf(\kcirc[T_1\. T_n][[T_{n+1}\. T_m]])$$ one defines $\gtX_\eta\subseteq\calM(k\{T_1\. T_m\})$ as the unit polydisc given by the conditions $|T_i|<1$ for $n<i\le m$. In particular, the polydisc is open when $n=0$ and is closed when $n=m$. For an affine $\gtY$ one defines $\gtY_\eta$ using a closed embedding into $\gtX$ as above, and for a general special formal scheme the functor is defined using gluing. An anti-continuous reduction map $\gtX_\eta\to\gtX_s$ is defined as earlier.

\begin{exer}
(i) $\gtX_\eta$ is a good $k$-analytic space, and it is strict when the valuation is non-trivial.

(ii) In the case of affine $\gtX=\Spf(A)$, the generic fiber can be identified with the set of real semivaluations on $A$ that extend the valuation on $k$, are bounded and are strictly smaller than one on the elements of an ideal of definition of $A$.

(iii) In the affine case, $\gtX_\eta$ is an increasing union of affinoid domains $X_n$ such that $X_n$ is Weierstrass in each $X_m$ for $m\ge n$ (e.g. an open polydisc is an increasing union of smaller closed polydiscs).
\end{exer}

One of motivations to introduce generic fibers of special fibers is the following result of Berkovich.

\begin{fact}\label{formber}
Let $\gtX$ be a $\kcirc$-special formal scheme (e.g. a formal scheme of finite type over $\kcirc$) and let $Z\into\gtX_s$ be a closed subscheme. Then the preimage of $Z$ under the reduction map $\pi_\gtX{\colon}\gtX_\eta\to\gtX_s$ depends only on the formal completion of $\gtX$ along $Z$. Moreover, this preimage is precisely the generic fiber $(\hatgtX_Z)_\eta$ of the formal completion of $\gtX$ along $Z$.
\end{fact}

The following conjecture in the opposite direction describes the precise information about the formal scheme that is kept in the generic fiber $(\hatgtX_Z)_\eta$.

\begin{conj}\label{formconj}
If $\gtX$ is locally of the form $\Spf(\calAcirc)$ for $k$-affinoid algebras $\calA$ (i.e. $\gtX$ is of finite type over $\kcirc$ and is normal in its generic fiber), then the henselization of $\gtX$ along a closed subscheme $Z\into\gtX_s$ is completely determined by the generic fiber $(\hatgtX_Z)_\eta$.
\end{conj}

A partial evidence in favor of this conjecture is provided by the following result that was proved in \cite{insepunif} and applied to resolution of singularities in positive characteristic.

\begin{fact}\label{formtem}
Assume that $\gtX'$ and $\gtX$ are locally of the form $\Spf(\calAcirc)$ for $k$-affinoid algebras $\calA$, $Z\into\gtX_s$ and $Z'\into\gtX'_s$  are closed subschemes, and $f:\gtX'\to\gtX$ is a morphism that induces isomorphism between the preimages of $Z'$ and $Z$ in the generic fibers. Then the restriction of $f$ onto a small neighborhood of $Z'$ is strictly \'etale over $Z$.
\end{fact}

\subsection{Raynaud's theory}\label{raysec}

\subsubsection{An overview}
Assume that the valuation is non-trivial. We constructed a functor $\eta$ whose source is the category $\kcirc$-$Fsch$ of formal $\kcirc$-schemes of finite type and whose target is the category $st$-$k$-$An^c$ of compact strictly $k$-analytic spaces. Raynaud's theory completely describes this functor in the following terms: $\eta$ is the localization of the source by an explicitly given family of morphisms $\calB$. In particular, one can view a compact strictly $k$-analytic space $\gtX_\eta$ as its {\em formal model} $\gtX$ given up to a morphism from $\calB$. (A possible analogy is to think about $\gtX$ as a particular atlas of a manifold $\gtX_\eta$ with morphisms from $\calB$ being the refinements of the atlases.)

Clearly, the central part of the theory should be to describe the morphisms $f{\colon}\gtY\to\gtX$ that are rig-isomorphisms (or generic isomorphisms). Although this is not so easy, we will find a nice cofinal family $\calB$ among all rig-isomorphisms. To guess what such $\calB$ can be, let us consider a very similar problem in the theory of schemes. Given a scheme $X$ with a schematically dense open subscheme $U$ (which will play the role of the generic fiber) by a $U$-modification of $X$ we mean a proper morphism $f{\colon}X'\to X$ such that $f^{-1}(U)$ is schematically dense in $X'$ and is mapped isomorphically onto $U$. A strong version of Chow lemma states that the family of $U$-admissible blow ups, i.e. blow ups whose center is disjoint from $U$, form a cofinal family among all $U$-modifications. Now, it is natural to expect that in our situation one can take $\calB$ to be the family of all formal blow ups along open ideals (i.e. ideals supported on $\gtX_s$).

\subsubsection{Admissible blow ups}

\begin{defin}
Recall that the blow up $\Bl_\calJ(X)$ of a scheme $X$ along an ideal $\calJ\subseteq\calO_X$ is defined as $\bfProj(\oplus_{n=0}^\infty\calJ^n)$. If $A$ is an $I$-adic ring then the formal blow up of the affine formal scheme $\Spf(A)$ along an ideal $J\subseteq A$ is defined as the $I$-adic completion of $\Bl_J(\Spec(A))$. This definition is local on the base and hence globalizes to a definition of {\em formal blow up} $\hatBl_\calJ(\gtX)$ of a formal scheme $\gtX$ along an ideal $\calJ\subseteq\calO_\gtX$. If $\calJ$ is open (i.e. contains an ideal of definition) then the formal blow up is called {\em admissible}.
\end{defin}

\begin{fact}
(i) Any composition of (admissible) formal blow ups is an (admissible) formal blow up.

(ii) (Admissible) blow ups form a filtered family.
\end{fact}

\begin{exer}
If $\gtX$ is of finite type over $\kcirc$ then any admissible formal blow up is a rig-isomorphism.
\end{exer}

\subsubsection{The main results}

\begin{fact}[Raynaud]\label{Rfact}
The family $\calB$ of formal blow ups in the category $\kcirc$-$Fsch$ admits a calculus of right fractions and the localized category is equivalent to $st$-$k$-$An^c$. The localization functor is isomorphic to the generic fiber functor.
\end{fact}

\begin{rem}
This fact implies the following two corollaries:

(i) The family of admissible blow ups of a formal scheme $\gtX$ from $\kcirc$-$\Fsch$ is cofinal in the family of all rig-isomorphisms $\gtX'\to\gtX$.

(ii) Each compact strictly $k$-analytic space $X$ admits a {\em formal model} $\gtX$, i.e. a formal scheme $\gtX$ in $\kcirc$-$\Fsch$ with an isomorphism $\gtX_\eta\toisom X$.

Actually, these two statements serve as intermediate steps while proving Fact \ref{Rfact}. Moreover, one proves that if $\{X_i\}$ is a finite family of compact strictly analytic domains in $X$ then there exists a model $\gtX$ with open subschemes $\gtX_i$ such that $(\gtX_i)_\eta\toisom X_i$.
\end{rem}

Let us say a couple of words on the proof of Fact \ref{Rfact}. The functor $\eta$ takes morphisms of $\calB$ to isomorphisms hence it induces a functor $\calF{\colon}\kcirc$-$Fsch/\calB\to st$-$k$-$An^c$. One easily sees that $\calF$ is faithful. The proof that $\calF$ is full reduces to proving claim (i) of the above remark. This is essentially a strong version of the Chow lemma and the main ingredient in its proof is Gerritzen-Grauert theorem. Once we know that $\calF$ is fully faithful, it remains to show that it is essentially surjective, i.e. to prove claim (ii) of the remark. For a strictly $k$-affinoid space it is very easy to find a formal model, and in general one chooses an affinoid covering $X=\cup_i X_i$, finds models $\gtX_i$ and then uses that $\calF$ is full to find formal blow ups $\gtX'_i\to\gtX$ so that the formal schemes $\gtX'_i$ glue to a formal model $\gtX'$ of $X$.

\subsection{Rigid geometry}
A naive attempt to construct the generic fiber of an affine formal scheme $\gtX=\Spf(A)$ is to declare that $\gtX_\eta$ is the set $\Spec(A)\setminus\Spf(A)$ of all non-open ideals. Such definition is not compatible with formal localizations, because the Zariski topology becomes stronger on localizations. In particular, this definition cannot be globalized. The situation improves, however, if one only considers the set of closed points of $\Spec(A)\setminus\Spf(A)$. In some cases such spectrum can be globalized, though the usual Zariski topology should be replaced with a certain $G$-topology in order to achieve this. We will not develop this point of view, but we will show how a similar approach gives rise to the rigid geometry of Tate.

Let $k$ be a non-archimedean field with a non-trivial valuation. Strictly $k$-affinoid algebras are the basic objects of rigid geometry over $k$. An affinoid space $X_0=\Sp(\calA)$ is defined as the set of maximal ideals of $A$ provided with the $G$-topology of finite unions of affinoid domains. By Hilbert Nullstellensatz, the residue field of any point $x\in X_0$ is finite over $k$ and hence $X_0$ is the set of Zariski closed points of $X=\calM(\calA)$. The theory of rigid affinoid spaces and general rigid analytic spaces is developed similarly to the theory of strictly $k$-analytic spaces from \S\S{2--3}. Some intermediate results are slightly easier to
prove because we only worry for Zariski closed points, but in the end one has less tools to solve problems. For example, Shilov boundaries and the class of good spaces are not seen in rigid geometry. Another example of an application where generic points are very important is the theory of \'etale cohomology of analytic spaces. There exist non-zero \'etale sheaves which have zero stalks at all rigid points but they necessarily have a non-zero stalk at a point of $X$. (Rigid points form a conservative family for coherent sheaves, and so the latter are easily tractable in the framework of Tate's rigid geometry.)

\subsection{Adic geometry}\label{adicsec}
\subsubsection{} Adic geometry replaces formal schemes with more general objects that have a honest generic fiber (as an adic space). Let us recall why formal schemes have no generic fiber. Let $k$ be a non-archimedean field with a non-trivial valuation and non-zero $\pi\in\kcirccirc$ and let $A$ be a $\pi$-adic $\kcirc$-algebra of finite type over. Formal inverting of $\pi$ produces the zero ring in two stages: first we invert $\pi$ obtaining a strictly $k$-affinoid algebra $\calA$ and then we have to factor over the unit ideal because the $\pi$-adic topology on $\calA$ is trivial (and so the $\pi$-adic separated completion of $\calA$ is $0$). This suggests to extend the category of adic rings so that topological rings like $\calA$ (with its Banach topology) are included. R. Huber suggested a way to do that, and it is very natural if we recall how the topology of $k$ is actually defined.

\begin{defin}
An {\em $f$-adic ring} is a topological ring that contains an open adic ring $A_0$ with a finitely generated ideal of definition. Any such $A_0$ is called a {\em ring of definition} (because it can be used to define the topology of $A$).
\end{defin}

Note that ring of definition is an analog of a unit ball for a norm.

\begin{exer}
(i) Any $k$-affinoid algebra $\calA$ is $f$-adic.

(ii) $\calA$ is reduced if and only if $\calAcirc$ is a ring of definition.
\end{exer}

Adic spectrum is defined analogously to analytic spectrum but using all continuous semivaluations. This forces one to modify the notion of a basic ring as follows.

\begin{defin}
(i) An {\em affinoid ring} is a pair $A=(A^\rhd,A^+)$ where $A^\rhd$ is an $f$-adic ring and $A^+$ is an open subring which is integrally closed in $A^\rhd$ and is contained in the ring of power-bounded elements $A^\circ$. The ring $A^+$ is called the {\em ring of integers} of $A$.

(ii) The adic spectrum $\Spa(A)$ is the set of all equivalence classes of continuous semivaluations on $A^\rhd$ such that $|a|\le 1$ for any $a\in A^+$.
\end{defin}

The (usual!) topology and the structure sheaf on $\Spa(A)$ are defined using rational domains. This is done similarly to our definitions, so we omit the details.

\begin{exex}
(i) The $k$-affinoid rings in adic geometry are the pairs $(\calA,\calAcirc)$ with $(k,\kcirc)$ playing the role of the ground field. Show that the only $\kcirc$-ring of integers of a strictly $k$-affinoid algebra $\calA$ is $\calAcirc$. (Hint: use Exercise \ref{torex}.)

(ii) For any adic $\kcirc$-algebra $A$ of finite type, $(A,A)$ is an affinoid ring. Check that $\gtX=\Spa(A,A)$ is an adic space and set-theoretically $\gtX$ is the disjoint union of the {\em generic fiber} $\gtX_\eta=\Spa(\calA,\calAcirc)$ where $\calA=A_\pi$ and the closed fiber $\gtX_s$ which consists of all semivaluations that vanish on $\pi$. Note that $\Spf(A)$ naturally embeds into $\gtX_s$ as the set of points of $\gtX$ with trivial valuation on the
residue field.

(iii) An {\em affinoid field} is an affinoid ring $k$ such that $k^\rhd$ is a valued field of height $h^\rhd\le 1$ (with the induced topology) and $k^+$ is a valuation ring of $k$ contained in $\kcirc$. Let $h^+$ be the height of $k^+$. Show that $\Spa(k)$ is a chain (under specialization) of $h^++1-h^\rhd$ points.
\end{exex}

\begin{rem}
(i) Spectra of affinoid fields are "atomic objects" in the sense that they do not admit non-trivial monomorphisms from other spaces such that the closed point is contained in the image. In particular, a point of height at least two or a fiber of a morphism over such point is not an adic space.

(ii) Points of height at least two are very different from the usual analytic points. For example, their local rings are usually not henselian (because there is no reasonable completion for valued fields of height more than one).
\end{rem}

\begin{exer}
Describe all adic points of the affine line over $k$. (Hint: show that the only new points are height two points contained in the closures of type two points $x\in\bfA^\an_k$. Each connected component of $\bfA^\an_k\setminus\{x\}$ contains one such point.)
\end{exer}

\subsection{Comparison of categories and spaces}
\subsubsection{} In principle, all approaches to non-archimedean analytic geometry produce the same category of spaces of finite type (with rational radii of convergence).

\begin{fact}
The following categories are naturally equivalent: (a) the category of compact strictly $k$-analytic spaces, (b) the category of formal $\kcirc$-schemes of finite type localized by admissible blow ups, (c) the category of quasi-compact and quasi-separated rigid $k$-analytic space, (d) the category of quasi-compact and quasi-separated adic $\Spa(k,\kcirc)$-spaces of locally finite type.
\end{fact}

Let $\gtX$ be a formal $\kcirc$-scheme of finite type and let $\gtX_\eta^\rig$, $\gtX_\eta^\an$ and $\gtX_\eta^\ad$ be its generic fibers in the three categories of non-archimedean spaces over $k$. On the level of topological spaces these objects are related as follows.

\begin{fact}
One has that $\projlim_{f{\colon}\gtX'\to\gtX}|\gtX'|\toisom|\gtX_\eta^\ad|\supset|\gtX^\an_\eta|\supset|\gtX^\rig_\eta|$,
where the limit is taken over all admissible formal blow ups $f$. Furthermore, $\gtX^\an_\eta$ is the set of all points of height one in $\gtX_\eta^\ad$ and also it is homeomorphic to the maximal Hausdorff quotient of $\gtX^\ad_\eta$, and $\gtX^\rig_\eta$ is the set of Zariski closed points of $\gtX^\an_\eta$.
\end{fact}

Finally, the sheaves on these spaces are connected as follows.

\begin{fact}
The topoi (i.e. the categories of sheaves of sets) of the following sites are equivalent: $\gtX^\an_\eta$ with the $G$-topology of compact strictly $k$-analytic domains, $\gtX^\rig_\eta$ with the topology of compact rigid domains, and $\gtX_\eta^\ad$ with its usual topology. In particular, $\gtX_\eta^\ad$ is simply the set of (equivalence classes of) points of all these sites.
\end{fact}

\subsection{Reduction of germs and Riemann-Zariski spaces}\label{RZsec}

The aim of this section is to use certain reduction data to study the local structure of an analytic space $X$ at a point $x$. For simplicity we will assume that $X$ is strictly analytic, and will discuss the general case in \S\ref{gensec}. The main tool we will develop and use is a  reduction functor that associates to the germ $(X,x)$ certain birational object $(X,x)^\sim$ over $\tilk$. We will see that it has interesting connections to Raynaud's theory, adic geometry and the classical (though not very well known) Riemann-Zariski spaces. The reduction functor is used to establish fine local properties of analytic spaces, including Facts \ref{boundfact} and \ref{faithfact}.

\subsubsection{Reduction of germs: preliminary version}

To motivate the definition of $(X,x)^\sim$ let us first pursue the following mini-goal: find a simple criterion for an analytic space $X$ to be not good at a point $x\in X$ and use it to justify the claim of Example \ref{nogoodexam}. Note that an affinoid space $X=\calM(A)$ possesses an affine formal model $\gtX$: although $\calAcirc$ dos not have to be of finite type over $\kcirc$, we can simply take an admissible surjection $f:k\{\uT\}\to\calA$ and set $\gtX=\Spf(f(\kcirc\{\uT\}))$. By Raynaud's theory, for any other formal model $\gtX'$ there exists a formal model $\gtX''$ which is an admissible blow up of both $\gtX$ and $\gtX'$. It follows that the reduction $\tilX=\gtX'_s$ is a $\tilk$-variety of a rather special form: there exists an $\tilX$-proper variety which admits a proper morphism to an affine variety. It turns out that this observation can be localized to a criterion of goodness at a point.

\begin{exer}\label{goodcrit}
(i)* Let $X$ be strictly analytic, let $x\in X$ be a point and let $\gtX$ be a formal model with reduction $\tilX=\gtX_s$. Let $Z$ denote the Zariski closure of the image of $x$ under the reduction map $\pi_\gtX:X\to\tilX$, and assume that there is no proper morphism $Z'\to Z$ such that $Z'$ admits a proper morphism to an affine variety. Show that $X$ is not good at $x$. (Hint: use Fact \ref{intboufact}(vi) and Raynaud's theory.)

(ii) Deduce that both spaces in Example \ref{nogoodexam} are non-good.
\end{exer}

Now, it is natural to expect that one can define an interesting reduction invariant of the germ $(X,x)$ by considering varieties $Z=\ol{\pi_\gtX(x)}$ as above modulo an equivalence relation that eliminates the freedom in the choice of the reduction. A straightforward attempt would be to factor the category of $\tilk$-varieties by surjective proper morphisms, but we can slightly refine this by taking the residue field $\wHx$ into account (recall that it does not have to be finitely generated over $\tilk$). This is based on the observation that $\pi_\gtX(x)$ is the image of the reduction morphism $\Spec(\wHx)\to\tilX$. Now we are ready to give a preliminary definition of the reduction functor.

\begin{defin}\label{birdef}
(i) Let $\Var_\tilk$ be the category of {\em pointed $\tilk$-varieties} (i.e. morphisms $\Spec(K)\to Z$, where $K$ is a $\tilk$-field and $Z$ is a $\tilk$-variety) with morphisms $f$ been compatible pairs of morphisms $f_\eta:\Spec(K')\to\Spec(K)$ and $f_s:Z'\to Z$. Let $\calB$ denote the family of morphisms in which $f_s$ is proper and $f_\eta$ is an isomorphism, and let $\Var_\tilk/\calB$ denote the localization category.

(ii) For any compact strictly analytic $X$ with a point $x$ and a formal model $\gtX$ consider the object $\Spec(\wHx)\to\gtX_s$ of $\Var_\tilk$. The same object viewed in the category $\Var_\tilk/\calB$ is denoted $(X,x)_\gtX^\sim$
\end{defin}

\begin{fact}\label{birfact}
The isomorphism class of $(X,x)_\gtX^\sim$ depends only on the germ $(X,x)$ (and so is independent of the choice of $\gtX$ and is preserved if we replace $X$ with a neighborhood of $x$). Moreover, this construction is functorial on the category of germs of strictly analytic spaces at a point.
\end{fact}

\subsubsection{The category $\bir_\tilk$}\label{birsec}
The family $\calB$ in Definition \ref{birdef} is easily seen to admit the calculus of right fractions (i.e. any morphism in $\Var_\tilk/\calB$ is of the form $f\circ b^{-1}$ where $b\in\calB$ and $f$ is a morphism in $\Var_\tilk$). In particular, the localization functor can be easily described. Nevertheless, it would be desirable to have a more geometric interpretation of $\Var_\tilk/\calB$, and it turns out that the latter is provided by classical RZ (or Riemann-Zariski) spaces.

\begin{defin}
For a $\tilk$-field $K$ let $\RZ_\tilk(K)$ denote the set of all valuation rings $k\subseteq\calO\subseteq K$ with $\Frac(\calO)=K$. Provide $\RZ_\tilk(K)$ with the topology whose basis is formed by the sets $$\RZ_\tilk(K[f_1\. f_n])=\{\calO\in\RZ_\tilk(K)|\ f_1\. f_n\in\calO\}$$ for any choice of $f_1\. f_n\in K$.
\end{defin}

\begin{rem}
Such spaces were introduced by Zariski in 1930ies. He called them Riemann spaces for their (very relative) analogy with Riemann surfaces. These spaces are sometimes used in birational geometry, for example, for resolution of singularities, and their modern name is Riemann-Zariski or Zariski-Riemann spaces (in both variants).
\end{rem}

\begin{fact}
The spaces $\RZ_\tilk(K)$ are qcqs (i.e. quasi-compact and quasi-separated).
\end{fact}

\begin{rem}
This simple fact is due to Zariski. Funny enough, it was reproved in many works and often incorrectly (including Nagata's work on compactification and \cite{temred1}). The mistake is always the same -- one assumes that quasi-compactness is preserved under projective limits, which is incorrect in general.
\end{rem}

\begin{defin}
A {\em birational space over $\tilk$} is a qcqs topological space $X$ with a local homeomorphism $\phi:X\to\RZ_\tilk(K)$. A morphism between $\phi':X'\to\RZ_\tilk(K')$ and $\phi$ is a $\tilk$-embedding $i:K\into K'$ and a continuous map $X'\to X$ compatible with the map $\RZ_\tilk(K')\to\RZ_\tilk(K)$ induced by $i$. The category of birational spaces over $\tilk$ will be denoted $\bir_\tilk$.
\end{defin}

\begin{rem}
In this definition, $\phi$ plays the role of a structure sheaf. Such a naive structure sheaf suffices to define the category $\bir_\tilk$. For the sake of completeness we note that one can provide $X$ with a real structure sheaf $\calO_X$ whose stalk at a point $x\in X$ is the corresponding valuation ring. This more refined approach becomes useful when one works with (much more general) relative Riemann-Zariski spaces.
\end{rem}

\begin{fact}\label{smallray}
For any pointed $\tilk$-variety $f:\Spec(K)\to X$ of $\Var_\tilk$ one can construct a birational space consisting of all pairs $(\calO,g)$, where $\calO\in\RZ_\tilk(K)$ and $g:\Spec(\calO)\to X$ is a morphism compatible with $f$ (forgetting $g$ defines an obvious projection onto $\RZ_\tilk(K)$). This defines a functor $\Val:\Var_\tilk\to\bir_\tilk$ which induces an equivalence $\Var_\tilk/\calB\toisom\bir_\tilk$.
\end{fact}

\begin{defex}
(i) Given an object (resp. morphism) $X$ of $\bir_\tilk$, by a {\em scheme model} of $X$ we mean any object $\calX$ (resp. morphism) of $\Var_\tilk$ such that $\Val(\calX)\toisom X$.

(ii) A birational space $\phi:X\to\RZ_\tilk(K)$ is {\em affine} if $\phi$ is injective and its image is of the form $\RZ_\tilk(K[f_1\. f_n])$. Show that $\phi$ is affine if and only if it admits an affine scheme model.

(iii) A morphism $g:\phi'\to\phi$ is {\em separated} (resp. {\em proper}) if the map $X'\to X\times_{\RZ_\tilk(K)}\RZ_\tilk(K')$ is injective (resp. bijective). Show that this happens if and only if $g$ admits a separated (reps. proper) scheme model, and then any scheme model is separated (resp. proper). (Hint: this is a slightly refined version of the valuative criteria.)

(iv) A birational space $\phi:X\to\RZ_\tilk(K)$ is {\em separated} (resp. {\em proper}) if so is its morphism to the final object $\RZ_\tilk(\tilk)$. Show that this happens if and only if $\phi$ is injective (resp. bijective).

(v)* Show that there is no reasonable notion of affine morphisms in $\bir_\tilk$. (Hint: if you have solved Exercise \ref{goodcrit} then you may already know an appropriate example of a morphism in $\bir_\tilk$.)
\end{defex}

The following fact is already due to Zariski.

\begin{fact}
If $X\to\RZ_\tilk(K)$ is a birational space and $\Spec(K)\to\calX_i$, $i\in I$ is the family of all its scheme models then the natural map $X\to\projlim_{i\in I}\calX_i$ is a homeomorphism.
\end{fact}

Finally, we define germ reductions as birational spaces.

\begin{defin}
For any germ $(X,x)$ of a strictly analytic space $X$ take a compact neighborhood $V$ of $x$ with a formal model $\gtV$ and define $(X,x)^\sim$ to be the image of $(V,x)^\sim_\gtV$ in $\bir_\tilk$. It is called {\em reduction of $X$ at $x$} or {\em reduction of the germ $(X,x)$}.
\end{defin}

Due to Facts \ref{birfact} and \ref{smallray}, germ reduction is a functor from the category of germs to the category of birational spaces.

\subsubsection{Relation to other theories}

\begin{rem}
There is a strong analogy between formal $\kcirc$-schemes, adic spaces and the functor $\gtX\mapsto\gtX_\eta^\ad$ on one side, and pointed varieties, birational spaces and the functor $\Val$ on the other side. (Also, we stress this analogy by saying formal models and scheme models.) In particular:

(a)Both adic and birational spaces are built from valuations of arbitrary height.

(b) Raynaud's theory says that $\eta$ induces an equivalence of a category of formal schemes localized by admissible blow ups and a category of adic spaces. The equivalence $\Var_\tilk/\calB\toisom\bir_\tilk$ induced by $\Val$ can be viewed as a baby version of Raynaud's theory.

(c) On the set-theoretic level, Raynaud's theory reduces to homeomorphism between adic (resp. birational) spaces and projective limits of their formal (resp. scheme) models.
\end{rem}

An additional link between adic spaces and the germ reduction functor is provided by the following fact.

\begin{fact}\label{redfact}
Let $X$ be a strictly analytic space with associated adic space $X^\ad$, let $x\in X$ be a point and let $\ox$ be the closure of $x$ in $X^\ad$. Then $(X,x)^\sim$ is naturally homeomorphic to $\ox$.
\end{fact}

\begin{rem}
(i) This fact suggests an alternative definition of the reduction functor. It is more elegant, but not suited for computations.

(ii) In some sense, the role of the reduction functor is to make visible some local information contained in $X^\ad$ but not seen directly in the analytic geometry.
\end{rem}

\subsubsection{Main properties and applications}
It turns out that many local properties of analytic spaces and their morphisms are reflected by the reduction functor, and here is the list of the main ones.

\begin{fact}\label{mainredfact}
(i) The reduction functor establishes a bijection between germ subdomains $(Y,x)$ of $(X,x)$ and birational subspaces of $(X,x)^\sim$. This bijection preserves intersections, finite unions and inclusions.

(ii) A strictly analytic space $X$ is good at a point $x$ if and only if the reduction $(X,x)^\sim$ is affine.

(iii) A morphism $f:Y\to X$ of strictly analytic spaces is separated at a point $y\in Y$ (resp. $y\in\Int(Y/X)$) if and only if the reduction morphism $(Y,y)^\sim\to(X,f(y))^\sim$ is separated (resp. proper).
\end{fact}

\begin{exer}
Check that Fact \ref{mainredfact}(iii) reduces to Fact \ref{intboufact}(vi) when $X$ and $Y$ are affinoid. In particular, show that if $X=\calM(\calA)$ and $x\in X$ is a point with character $\chi_x:\calH(x)\to\calA$ then $\Spec(\wHx)\to\Spec(\tilchi_x(\tilcalA))$ is a scheme model of $(X,x)^\sim$.
\end{exer}

Note that Fact \ref{mainredfact}(ii) gives a simple necessary and sufficient criterion of goodness. We have already used this criterion to show that certain spaces are not good. The opposite implication (i.e. the criterion of goodness) is the difficult one, and it is the deepest property asserted by the Fact. All other claims follow relatively easily. (We already mentioned in Remark \ref{proprem} that the most difficult task in studying properness is to show that certain spaces are good.)

\subsubsection{Germ reduction of non-strict spaces}\label{gensec}
Finally, let us discuss how the reduction functor can be extended to non-strict analytic spaces. There exists no generalization of formal models and Raynaud's theory for general $H$-strict $k$-analytic spaces (though affinoid reduction can be defined as in Remark \ref{grrem}). On the other hand, one can define a category $\bir_{\tilk_H}$ of $H$-graded birational spaces over $\tilk_H$ analogous to the category $\bir_\tilk$ (this requires to define $H$-graded valuation rings, etc.). Then an $H$-graded germ reduction functor with values in $\bir_{\tilk_H}$ can be constructed as follows: for an affinoid $X=\calM(\calA)$ take $(X,x)_H^\sim$ to be the $H$-graded birational space corresponding to the homomorphism $\wt{(\chi_x)}_H:\tilcalA_H\to\wHx_H$, and in general we cover a germ $(X,x)$ by good subdomains $(X_i,x)$ and glue $(X,x)_H^\sim$ from $(X_i,x)^\sim_H$.

\begin{fact}\label{f}
(i) All assertions of Fact \ref{mainredfact} generalize verbatim to the $H$-graded setting.

(ii) There is a natural fully faithful embedding $\bir_{\tilk_H}\into\bir_{\tilk_G}$, where $G=\bfR^\times_+$. An analytic space $X$ is $H$-strict locally at $x$ if and only if the reduction $(X,x)^\sim_G$ comes from the subcategory $\bir_{\tilk_H}$.
\end{fact}

Fact \ref{f}(ii) provides a local criterion of $H$-strictness which was proved in \cite{temred2} for strictly analytic category and generalized to any $H$ in \cite{ct}. As a simple corollary, it was shown in loc.cit. that $k_H$-$An$ is a full subcategory of $k$-$An$.

\section{Analytic curves}\label{cursec}\label{6}
\begin{definsect}
(i) A {\em $k$-analytic curve} $C$ is a $k$-analytic space of pure dimension one, i.e. $\dim(C)=1$ and $C$ does not contain discrete Zariski closed points.

(ii) In the same way as in Definition \ref{typesex} we divide the points of $C$ into four types accordingly to their completed residue field.
\end{definsect}

\subsection{Examples}

\subsubsection{Constructions}
Let us first list some constructions that allow to create/enrich our list of $k$-analytic curves: (i) analytification of an algebraic curve, (ii) generic fibers of formal curves, (iii) an analytic domain in a curve, (iv) a finite covering of a curve (or, more generally, a covering with discrete fibers).

\begin{exex}
(i) The following curves can be obtained by the first method: affine line, projective curves, affine line with doubled origin.

(ii) Let $\calX$ be an irreducible projective algebraic $k$-curve with $k(\calX)=K$. The Zariski closed points of $X=\calX^\an$ are in one-to-one correspondence with the closed points of $\calX$, and other points of $X$ are in one-to-one correspondence with the valuations on $K$ that extend the valuation on $k$. In particular, $K\into\calH(x)$ and $\hatK\toisom\calH(x)$ for any non Zariski closed point $x\in X$.

(iii) The following curves can be obtained by the second method: compact $k$-analytic curves.

(iv) Most of Hausdorff curves admit a formal model of locally finite type over $\kcirc$. For example, find such models for an affine line and for an open unit disc.
\end{exex}

Now, let us study the other two methods with more details.

\subsubsection{Domains in the affine line}
A typical example of a compact domain $X$ is a closed disc $E(a,r)$ with finitely many removed open discs $E(a_i,r_i)$.

\begin{exer}\label{irrex}
(i) Prove that $X$ is a Laurent domain in $E(a,r)$.

(ii) Show that if $r$ and $r_1$ are linearly independent over $|k^\times|$ then $X$ is not a finite covering of a disc. (Hint: if $X=\calM(A)$ is finite over $E(b,s)$ then
$\rho(\calA)\subset\{0\}\cup\sqrt{s^\bfZ|k^\times|}$.)

(iii) Show that one can extend $X$ a little bit so that $r_i\in\sqrt{r^\bfZ|k^\times|}$ and then $X$ is a finite covering of $E(0,r)$.
\end{exer}

Next we describe neighborhoods of an especially simple form.

\begin{exer}
Assume that $k=k^a$. Show that a point $x\in\bfA^1_k$ admits a fundamental family of open neighborhoods $X_i$ as follows:

(i) if $x$ is of type 1 or 4 then $X_i$ are open discs,

(ii) if $x$ is of type 3 then $X_i$ are open annuli,

(iii) if $x$ is of type 2 then $X_i$ are open discs with removed finitely many closed discs, and in addition one can achieve that $X_i\setminus\{x\}$ is a disjoint union of open discs and finitely many open annuli.
\end{exer}

Finally, let us discuss some open domains in $\bfA^1_k$ with $k=k^a$.

\begin{exer}
(i) Show that a separated gluing of two open annuli along an open annulus is an open annulus. (Hint: use Exercise \ref{glueex}(i).)

(ii) Show that a filtered union of a countable family of annuli does not have to be an annulus. (Hint: take $\bfP^1_k$ and remove two type 4 points.)

(iii) We say that a non-archimedean field $k$ is {\em local} if $|k^\times|\toisom\bfZ$ and $\tilk$ is finite. Show that either $k$ is finite over $\bfQ_p$ or $k\toisom\bfF_p((t))$. Drinfel'd upper half-plane is defined as $\bfP^1_k\setminus\bfP^1_k(k)$. Show that it is an open analytic domain in $\bfP^1_k$.
\end{exer}

\subsubsection{Finite covers}
First, we consider an inseparable cover giving rise to a pathology that cannot occur in the algebraic world.

\begin{exex}\label{badcov}
(i) Construct a non-archimedean field $k$ with $\cha(k)=p$ and $[k:k^p]=\infty$.

(ii) Choose elements $a_i\in k$ which are algebraically independent over $k$ and such that $|a_i|$ tend to zero. Set $\calA=k\{T,S\}/(S^p-\sum_{i=0}^\infty a_i T^{pi})$. Show that $X=\calM(\calA)$ is a finite covering of $E(0,1)$ of degree $p$, $X\otimes_k l$ is reduced for any finite field extension $l/k$ but $X\wtimes_k k^{1/p}$ is not reduced.
\end{exex}

Next, we study some quadratic covers. They will give rise to various interesting examples.

\begin{exer}
Assume that $\cha(\tilk)\neq 2$  and consider a series $f(T)=\sum_{i=0}^\infty a_iT^i\in k\{T\}$ with the affinoid algebra $\calA=k\{T,S\}/(S^2-f(T))$ and the quadratic covering $\phi{\colon}X=\calM(\calA)\to E(0,1)=\calM(k\{T\})$. Show that any type 4 point $x\in\bfA^1_k$ has two preimages, and the maximal point $p_r$ of the disc $E(0,r)$ has two preimages if and only if $|a_0|>|a_i|r^i$ for $i>0$. In particular, a Zariski closed point $x$ has two preimages if and only if $f(x)\neq 0$. (Hint: show that the binomial expansion of $\sqrt{1+z}$ has radius of convergence 1 over $k$.)
\end{exer}

Now let us study elliptic curves using double covers. For simplicity, we also assume that $k=k^a$.

\begin{exex}
It is known from algebraic geometry that any elliptic curve over $k$ can be realized as the double covering $\phi{\colon}E\to\bfP^1_k$ given by $S^2=T(T-1)(T-\lam)$.

(i) Assume that $|\lam|>1$.

(a) Show that the points with one preimage are precisely the points of the disjoint intervals $[0,1]$ and $[\lam,\infty]$. In particular, $X$ contains a cycle $\Delta(E)$ which is the preimage of the interval $I=[p_1,p_{|\lam|}]$ and a contraction of $\bfP^1_k$ onto $I$ lifts to the contraction of $E$ onto $\Delta(E)$.

(b) Show that the preimage of the disc $E(0,r)$ with $1<r<|\lam|$ is a closed annulus, and deduce that $E$ is glued from two annuli.

(c)* Show that $E$ can be obtained from $A(0;1,|\lam|^2)$ by identifying $A(0;1,1)$ with $A(0;|\lam|^2,|\lam|^2)$. Moreover, the universal cover of $E$ is isomorphic to $G_m$ and $G_m/q^\bfZ\toisom E$ for an element $q\in k$ with $|q|=|\lam|^2$.

(ii) Show that if $|\lam-1|<1$ or $|\lam|<1$ then the structure of $E$ is similar but with respect to other intervals connecting the four points.

(iii) Assume that $|\lam|=|\lam-1|=1$.

(a) Show that the points with one preimage are $p_1$ and the points of the disjoint intervals $[0,p_1)$, $[1,p_1)$, $[\lam,p_1)$ and $[\infty,p_1)$. Let $z=\Delta(E)$ be the preimage of $p_1$; then the contraction of $\bfP^1_k$ onto $p_1$ lifts to the contraction of $E$ onto $z$.

(b) Show that $E\setminus\{z\}$ is a disjoint union of open discs. Furthermore, $\wHz$ is of genus one over $\tilk$ and the closed points of its projective model parameterize the open discs of $E\setminus\{z\}$. In particular, $z$ is not locally embeddable into $\bfA^1_k$.
\end{exex}

The curves from (i) and (ii) are called Tate curves, or elliptic curves with bad reduction. The curves from (iii) are called curves with good reduction. In the sequel by {\em genus} of a type two point $z$ we mean the algebraic genus of $\wHz$ over $\tilk$. Points of positive genus are very special and very informative.

\begin{exer}
(i) Study curves $C$ of genus two given by $S^2=f(T)$ with $f(T)$ of degree five. Show that the first Betti number of $C$ plus the sum of genera of its type two points equals to the genus of $C$.

(ii)* Prove the same for any $C$ given by $S^2=f(T)$ where $f(T)$ is a polynomial without multiple roots.
\end{exer}

Finally, let us construct wild non-compact examples.

\begin{exer}\label{infex}
Assume that $k$ is not discretely valued.

(i) Show that if $|a_i|$ increase and tend to 1 then $f(T)=\sum_{i=0}^\infty a_iT^i$ is a bounded function with infinitely many roots on the open unit disc $D(0,1)$.

(ii) Show that by an appropriate choice of $f(T)$ as above one can achieve that the corresponding double cover $C$ of $D(0,1)$ has infinitely many loops and infinitely many positive genus points.
\end{exer}

It will follow from some further results that $C$ is an example of a non-compactifiable space. In the sequel we will study compactifiable (usually compact) curves.

\subsection{General facts about compact curves}

\subsubsection{Algebraization}

\begin{fact}
Any proper $k$-analytic curve $X$ is projective. In particular, $X$ is algebraizable.
\end{fact}

\begin{exer}
(i)* Prove the above fact. (Hint: take a Zariski closed point $P$ and show that $H^1(X,\calO_X(nP))$ vanishes for large enough $n$ by Kiehl's theorem on direct images \ref{kith}. Deduce that a linear system $\calO_X(nP)$ with large enough $n$ gives rise to a finite morphism from $X$ to the projectivization of $H^0(X,\calO_X(nP))$. Then use Fact \ref{gaga} from GAGA.

(ii) Deduce that the curve from Exercise \ref{badcov} is not a domain in a proper curve. Moreover, find $k$ with $[k:k^p]<\infty$ and a finite extension $K/K_r$ such that $K$ is not isomorphic to the completion of a finitely generated $k$-field of transcendence degree one. Use this to construct a $k$-affinoid curve that cannot be embedded into a proper curve.
\end{exer}

\subsubsection{Compactification}
We saw that if a separated compact curve $C$ is not geometrically reduced then it does not have to be embeddable into a proper one. In the opposite direction we have the following result.

\begin{fact}
Any separated geometrically reduced $k$-analytic curve is isomorphic to a domain in a projective $k$-curve.
\end{fact}

A general idea of the proof is as follows: we would like to patch the boundary of $C$, which consists of generic points (of types 2 and 3). One proves that a curve is geometrically reduced at a non Zariski closed point if and only if it is rig-smooth at such point. If $C$ is rig-smooth at $x$ then it admits a quasi-\'etale morphism $\phi{\colon}C\to\bfA^1_k$ locally around $x$. If we deform a quasi-\'etale morphism slightly then the isomorphism class of $C$ does not change. Therefore we can define $\phi$ using only equations of the form $\sum_{i=0}^n a_i(T)y^i$ where the coefficients $a_i$ are meromorphic. This allows to compactify $C$ at all
points of its boundary.

Using the Riemann-Roch theorem on a projective curve one deduces the following corollary.

\begin{fact}\label{affcur}
A separated, compact, and geometrically reduced $k$-analytic curve is affinoid if and only if it does not contain proper irreducible components.
\end{fact}

\subsubsection{Formal models}
In the sequel we assume that $C$ is a compact strictly $k$-analytic rig-smooth curve and the valuation is non-trivial. For any formal model $\gtC$ of $C$ let $\gtC^0\subset C$ be the preimage of the set of generic points of $\gtC_s$ under the reduction map.

\begin{fact}\label{formcur}
(i) The set $\gtC_0$ determines the formal model $\gtC$ up to a finite admissible blow up.

(ii) If $C$ is separated then a finite set $V$ of type 2 points is of the form $\gtC_0$ for some formal model if and only if $V$ contains the boundary of $C$ and hits each proper irreducible component of $C$.
\end{fact}

\begin{exer}
(i) Show that (ii) above does not hold in the non-separated case. (Hint: take the closed disc with doubled open disc, and patch in an open annulus instead of the doubled open disc. Then there is no formal model with a single generic point (although such a model exists as a not locally separated formal algebraic space).)

(ii)* Deduce Fact \ref{formcur} from Fact \ref{affcur}.
\end{exer}

\subsection{Rig-smooth curves}
Until the end of the paper, $C$ is a compact rig-smooth $k$-analytic curve. For simplicity, we also assume that $k$ is algebraically closed and $C$ is connected. In general, all our results hold up to a finite ground field extension and obvious corrections needed to deal with disconnected curves.

\subsubsection{Geometric structure of analytic curves}\label{geomsec}
Here is the main result about the structure of $C$ in its geometric formulation. An equivalent approach via formal models will be discussed later.

\begin{fact}\label{mainfact}
There exists a finite set $V$ of type two points such that $C\setminus V$ is a disjoint union of open discs and finitely many open annuli.
\end{fact}

This claim is very strong and implies many other important results that we state as exercises.

\begin{exer}
(i) $C$ has finitely many points of positive genus and $C$ can be contracted onto its subset $\Delta(V)$ homeomorphic to a finite graph. (Hint: take $\Delta(V)$ to be the the union of $V$ and the open intervals through the annuli; then $C\setminus\Delta(V)$ is a disjoint union of open unit discs which can be easily contracted.)

(ii)* If $C$ is proper then its algebraic genus equals to the sum of the genera of type two points plus the first Betti number of $\Delta(V)$. (Hint: use the semistable formal model associated to $V$ in the next section.)
\end{exer}

Next, let us study what is the freedom in the choice of $V$.

\begin{exer}\label{maincor}
(i) Cofinality: show that any finite set of type two points can be enlarged to a set $V$ as above.

(ii) Fact \ref{topconj} holds true for curves. (Hint: the sets $\Delta(V)$ form the required filtered family.)

(iii)* Minimality: show that there exists a minimal such $V$ unless $C\toisom\bfP^1_k$ or $C$ is a Tate curve. (Note that the degenerate cases are proper curves that can be covered by annuli.)
\end{exer}

Next, let us describe the local structure of $C$.

\begin{exer}
Show that a point $x\in C$ has a fundamental family of open neighborhoods $X_i$ such that

(i) $X_i$ are open discs when $x$ is of type 1 or 4,

(ii) $X_i$ are open annuli when $x$ is of type 2,

(iii) $X_i\setminus\{x\}$ are disjoint unions of open discs and finitely many open annuli when $x$ is of type 2.
\end{exer}

Using gluing of annuli and discs from Fact \ref{glueex}, it is easy to show that the above local description of $C$ is equivalent to its global description. This local fact, in its turn, easily reduces to study of the field $\calH(x)$. For example, for types 3 and 4 it suffices to show that $\calH(x)$ is topologically generated by an element $T$ (i.e. $\calH(x)=\widehat{k(T)}$). Surprisingly, no simple proof of this fact is known. A shortest currently known proof can be found in \cite{temst}, where it is used to obtain a new proof of the semistable reduction theorem (we will see that the semistable reduction theorem is equivalent to Fact \ref{mainfact}).

\subsubsection{Semistable formal models}

\begin{defin}
A formal $\kcirc$-scheme is {\em semistable} if it is \'etale-locally isomorphic to the formal schemes of the form $\gtZ_{n,a}=\Spf(\kcirc\{T_1\. T_n\}/(T_1\dots T_n-a))$ with $a\in\kcirc$.
\end{defin}

Let $\gtX$ be normal in its generic fiber. Then $\gtX$ is semistable if and only if it has the same formal fibers as the model schemes $\gtZ_{n,a}$ (this particular case of Conjecture \ref{formconj} is easily verified by a direct computation). In the case of curves this gives the following result.

\begin{fact}
A formal $\kcirc$-curve $\gtC$ with rig-smooth generic fiber is semistable if and only if the formal fibers over its closed points are open discs (over the smooth points) and open annuli (over the double points of $\gtC_s$).
\end{fact}

This exercise and Fact \ref{formcur} imply that the global description of $C$ given by Fact \ref{mainfact} is equivalent to the following fundamental result, which can be proved alternatively by a classical but rather complicated algebraic theory that involves stable reduction over a discretely valued field and the theory of moduli spaces of curves.

\begin{fact}[Semistable reduction theorem for analytic curves]
Any compact rig-smooth strictly analytic curve over an algebraically closed field $k$ admits a semistable formal model.
\end{fact}

In the same way, Exercise \ref{maincor} implies the following generalization of the above fact.

\begin{exer}
(i) Cofinality: any formal model $\gtC$ of $C$ admits an admissible blow up $\gtC'\to\gtC$ such that $\gtC'$ is semistable.

(ii) Stable reduction theorem: if $C$ is not isomorphic to $\bfP^1_k$ or a Tate curve then it possesses a minimal semistable formal model (called the stable formal model).
\end{exer}

\subsubsection{Skeletons}
The reader that solved Exercises in \S\ref{geomsec} is probably familiar with part of the ideas of this section. All facts of this section follow easily from the results of \S\ref{geomsec}.

\begin{defin}
(i) Let $V_0$ be a finite set of type 1 and 2 points of $C$. The {\em skeleton} $\Delta(C,V_0)$ is defined as follows: its set of vertices $V$ is the set of points $x\in C$ that are not contained in an open annulus $A\subset C\setminus V_0$, and its edges are formed by the points $x\in C$ that are not contained in an open disc $D\subset C\setminus V_0$. A vertex is {\em infinite} if it is of type 1.

(ii) The skeleton $\Delta=\Delta(C,V_0)$ is {\em degenerate} if the set $V_f$ of finite vertices is empty.
\end{defin}

For the sake of completeness, we make a remark about the more general situation that we do not study.

\begin{rem}
The definition makes sense for any curve $C$ over $k=k^a$ with a finite set $V$. For a curve over an arbitrary ground field $k$, its skeleton is defined as the image of the skeleton after the ground field extension to $\whka$.
\end{rem}

\begin{exer}
(i) Show that $\Delta$ is a finite graph whose infinite vertices are the points of $V_0$ of type 1.

(ii) Show that the only degenerate cases are when $C$ is a Tate curve and $V_0$ is empty, or $C\toisom\bfP^1_k$ and $V_0$ is a set of at most two type 1 points. Thus, in the degenerate cases $\Delta$ is empty, an infinite vertex, an interval with infinite ends, or a loop without vertices.

(iii) Show that in the non-degenerate case $V_f$ is the minimal set of points such that $C\setminus V_f$ is a disjoint union of open discs and annuli such that annuli are disjoint from $V_0$ and each open disc contains at most one infinite point of $V_0$.

(iv)* Show that if $C$ is not proper or $V_0$ has finite vertices then the above description of $\Delta$ implies (and is equivalent to) the following stable modification theorem: if $\gtC$ is a formal rig-smooth $\kcirc$-curve with a generically reduced Cartier divisor $\gtD$ then there exists a minimal {\em modification} $\gtC'\to\gtC$ (i.e. a proper morphism whose generic fiber is an isomorphism) such that $\gtC'$ is semistable and the strict transform of $\gtD$ is \'etale over $\Spf(\kcirc)$.

(v)* Formulate and prove an analogous statement when $C$ is proper and $V_0$ has no finite vertices. (Hint: this is the stable reduction theorem for a formal curve with a divisor.)
\end{exer}

\end{document}